\documentclass[11pt]{article}

\usepackage[utf8]{inputenc}
\usepackage{amssymb,amsmath}
\usepackage{amsmath,amscd}
\usepackage{amsthm}
\usepackage{graphicx}
\usepackage{amsrefs}
\usepackage{enumerate}
\usepackage{hyperref}

\usepackage{amsmath}
\usepackage{latexsym}
\usepackage{amsfonts}
\usepackage{amssymb}
\usepackage{makeidx}

\newcommand{\nf}{\hspace*{-5pt}}

\def \beq { \begin{equation} }
\def \eeq { \end{equation} }

\def\overset#1\to#2{\mathrel{\mathop{#2}\limits^{#1}}}
\def\underset#1\to#2{\mathrel{\mathop{#2}\limits_{#1}}}

\def \varinjlim {\underset{\longrightarrow}\to{lim}}
\def \varprojlim {\underset{\longleftarrow}\to{lim}}

\def \li {\varinjlim}
\def \lp {\varprojlim}

\def \lli#1 {\mathrel{\mathop{\li}\limits_{#1}}}
\def \llp#1 {\mathrel{\mathop{\lp}\limits_{#1}}}

\def \empty {\emptyset}
\def \leqs {\leqslant}
\def \rest {\restriction}

\def \im {\Rightarrow}

\def \hookr {\hookrightarrow }

\newcommand{\hfl}{\hspace*{\fill}}

\newlength{\dede}     

\newcommand{\nhp}{\hspace*{-\parindent}}

\newcommand{\N}{\mbox{$ \mathbb N$}}   

\newcommand{\lra}{\mbox{$\longrightarrow$}}  

\newcommand{\sub}{\mbox{$\subseteq$}}

\def \empty {\emptyset}

\newcommand{\cC}{\mbox{$\cal C$}}

\newcommand{\cI}{\mbox{$\cal I$}}

\newcommand{\cL}{\mbox{$\cal L$}}

\newcommand{\cS}{\mbox{$\cal S$}}
\newcommand{\cT}{\mbox{$\cal T$}}

\newcommand{\all}{\mbox{$\forall$}}

\newcommand{\Ra}{\mbox{$\Rightarrow$}}
\newcommand{\Lra}{\mbox{$\Leftrightarrow$}}

\def \rest {\restriction}
\def \ra {\rightarrow}

\def \im {\Rightarrow}

\def \hookr {\hookrightarrow }

\def\overset#1\to#2{\mathrel{\mathop{#2}\limits^{#1}}}
\def\underset#1\to#2{\mathrel{\mathop{#2}\limits_{#1}}}

\def \leqs {\leqslant}

\def \ll {\langle}
\def \rr {\rangle}

\def \varinjlim {\underset{\longrightarrow}\to{lim}}
\def \varprojlim {\underset{\longleftarrow}\to{lim}}

\def \li {\varinjlim}
\def \lp {\varprojlim}

\def \lli#1 {\mathrel{\mathop{\li}\limits_{#1}}}
\def \llp#1 {\mathrel{\mathop{\lp}\limits_{#1}}}

\newtheorem{Th}{Theorem}
\newtheorem{Co}[Th]{Corollary}
\newtheorem{Df}[Th]{Definition}
\newtheorem{Pro}[Th]{Proposition}
\newtheorem{Le}[Th]{Lemma}
\newtheorem{Exa}[Th]{Example}
\newtheorem{Rem}[Th]{Remark}
\newtheorem{Fa}[Th]{Fact}
\newtheorem{Que}[Th]{Question}
\newtheorem{Ct}[Th]{}
\newtheorem{Afi}[Th]{Claim}

\newcommand{\baf}{\begin{Afi}\nf{\sl }}
\newcommand{\eaf}{\end{Afi}}

\newcommand{\bdf}{\begin{Df}\nf{\bf }}
\newcommand{\edf}{\end{Df}}
\newcommand{\bte}{\begin{Th}\nf{\bf }}
\newcommand{\ete}{\end{Th}}  
\newcommand{\bco}{\begin{Co}\nf{\bf }}
\newcommand{\eco}{\end{Co}}
\newcommand{\ble}{\begin{Le}\nf{\bf }}
\newcommand{\ele}{\end{Le}}
\newcommand{\bpr}{\begin{Pro}\nf{\bf }}
\newcommand{\epr}{\end{Pro}}
\newcommand{\bex}{\begin{Exa}\nf{\bf } \rm}
\newcommand{\eex}{\end{Exa}}
\newcommand{\bre}{\begin{Rem}\nf{\bf } \rm}
\newcommand{\ere}{\end{Rem}}
\newcommand{\bfa}{\begin{Fa}\nf{\bf } \sl}
\newcommand{\efa}{\end{Fa}}
\newcommand{\bqt}{\begin{Que}\nf{\bf }}
\newcommand{\eqt}{\end{Que}}
\newcommand{\bdm}{\nhp{\bf Proof.  }}
\newcommand{\edm}{\smallskip}
\newcommand{\qdr}{\hfl $\square$ }
\newcommand{\bct}{\begin{Ct}\nf{\bf } \rm}
\newcommand{\ect}{\end{Ct}}
\newcommand{\bxa}{\begin{Exa}\nf{\bf } \rm}
\newcommand{\exa}{\end{Exa}}

\newcommand{\bu}{$\bullet$\ }

\newcommand{\no}{\mathit{No}}
\newcommand{\on}{\mathit{On}}

\newcommand{\pa}[2]{\langle #1, #2\rangle}

\DeclareMathOperator{\dom}{dom}

\newcommand{\sqsub}{\sqsubseteq}

\begin{document}

\title{An algebraic (set) theory of surreal numbers, I}
\author{Dimi Rocha Rangel \thanks{
Instituto de Matem\'atica e Estat\'istica, University of S\~ao Paulo, Brazil.
Emails: dimi@ime.usp.br, hugomar@ime.usp.br} \   Hugo Luiz Mariano}
\date{}
\maketitle

\begin{abstract}

The notion of surreal number was introduced by J.H. Conway in the mid
1970's: the surreal numbers constitute a linearly ordered (proper)
class $No$ containing the class of all ordinal numbers ($On$) that, 
working within the background set theory NBG, can be defined by
a recursion on the class $On$. Since then, have appeared
many constructions of this class and was isolated a full axiomatization of this
notion that been subject of interest due to large number of interesting properties
they have, including model-theoretic ones. Such constructions suggests strong
connections between the class $No$ of surreal numbers and the classes
of all sets and all ordinal numbers.

In an attempt to codify the universe of sets directly within the surreal number
class, we have founded some clues that suggest that this class is not suitable
for this purpose. The present  work, that expounds parts of the PhD thesis of the first author (\cite{Ran18}),  establishes a basis to obtain an "algebraic (set) theory
for surreal numbers" along the lines of the Algebraic Set Theory - a categorial set theory introduced in the 1990's based on the concept of ZF-algebra: 
to establish abstract and general links between the class of all
surreal numbers and a universe of "surreal sets" similar to the relations 
between the class of all ordinals ($On$) and the class of all sets ($V$), 
that also respects and expands the links between the linearly ordered class of all
ordinals and of all surreal numbers.

In the present work we introduce the notion of (partial) surreal algebra (SUR-algebra) 
and we explore some of its category theoretic properties, including (relatively) 
free SUR-algebras ($SA, ST$).

In a continuation of this work (\cite{RM19})
we will establish links, in both directions, 
between SUR-algebras and ZF-algebras (the keystone of  Algebraic Set Theory) and develop the first steps of  a certain kind of set theory based (or ranked)  on surreal numbers, that expands the relation between $V$ and $On$.

\end{abstract}

{\bf Keywords:} Surreal numbers; Algebraic Set Theory; SUR-algebras.

\section*{Introduction}

The notion of surreal number was introduced by J.H. Conway in the mid
1970's: the surreal numbers constitute a linearly ordered (proper)
class $No$ containing the class of all ordinal numbers ($On$) that, 
working within the background set theory NBG, can be defined by
a recursion on the class $On$. Since then, have appeared
many constructions of this class and was isolated a full axiomatization of this
notion. 

Surreal numbers have been subject of interest in many areas of Mathematics 
due to large number of interesting properties that they have:\\
- In  Algebra, through the concept of field of Hahn series and variants 
(see for instance \cite{mac39}, \cite{all62}, \cite{sct67}, \cite{DE01}, \cite{KM15});\\
- In Analysis (see the book \cite{Allbk});\\
- In Foundations of Mathematics, particularly in Model Theory,   since the surreal number
line is for proper class linear orders what the rational number line $\mathbb{Q}$ is for the countable: 
surreal numbers are the (unique up to isomorphism) proper class Fra\"iss\'e limit of the finite linearly ordered sets, they are set-homogeneous and universal for all proper class linear orders.

The plethora of aspects and applications of the surreals maintain the subject as an active research field. To make a point,  the 2016 edition of the
"Joint Mathematics Meetings AMS" --the largest Math. meeting in the world-- have counted 14 talks in its
"AMS-ASL Special Session on Surreal Numbers". \\
\url{http://jointmathematicsmeetings.org/meetings/national/jmm2016/2181_program_ss16.html}

Here we try to develop, from scratch, a new (we hope!) and  complementary foundational aspect of the Surreal Number Theory: to establish, in some sense,  a set theory based on the class of surreal numbers.

Set/class theories are one of the few fundamental mathematical theories that holds the power to base other notions of mathematics (such as points, lines, and real numbers). This is basically due to two aspects of these theories: the first is that the basic entities and relations are very simple in nature, relying only on the primitive notions of  set/class and a (binary) membership relation (``$ X \in Y $"), the second aspect is the possibility that this theory can perform  constructions of sets by several methods. This combination of factors allows to achieve a high degree of flexibility, in such a way that virtually all mathematical objects can be realized as being of some kind of set/class, and it has the potential to define arrows (category) as entities of the theory. In particular, the  set/class theories traditionally puts as a principle (the Axiom of Infinity) the existence of an "infinite set" - the smallest of these would be the set of all natural numbers - thus, such numbers are a derived (or "a posteriori) notion, which encodes the essence of the notion of "to be finite", that is apparently more intuitive.

The usual set/class theories (as ZFC or NBG) have the power of  "encode" (syntactically) its Model Theory:  constructions of models of set theory by the Cohen forcing method or through boolean valued models method are done by a convenient encoding of the fundamental binary relations $\in$ and $=$.

Let us list below some other fundamental theories:

\bu Set theories with additional predicates for non-Standard Analysis, as the E. Nelson's set theory named IST.

\bu P. Aczel's "Non-well-founded sets" (\cite{Aczbk}), where sets and proper classes are replaced by directed graphs (i.e., a class of vertices endowed with a binary relation)

\bu K. Lopez-Escobar "Second Order Propositional Calculus" (\cite{le09}), a theory with 
three primitive terms, that  encodes the full Second Order Intuitionistic Propositional Calculus  also includes Impredicative Set Theory.

\bu Toposes, a notion isolated in the 1970's by W. Lawvere and M. Tierney, provide generalized set theories, from  the category-theoretic point of view.

\bu Algebraic Set Theory (AST), another categorial approach to set/class theory, introduced in the 1990's by A. Joyal and I. Moerdijk (\cite{JMbk}).\\
Algebraic Set Theory replaces the traditional use of axioms on pertinence by categorial relations, proposing the general study of "abstract class categories" endowed with a notion of "small fibers maps". In the same way that the notion of `` being finite '' is given a posteriori in ZFC, after guaranteeing an achievement of the Peano axioms - which axiomatizes the algebraic notion of free monoid in 1 generator - the notions of "to be a set" and  "be an ordinal" are given a posteriori in AST. The class of all sets
is determined by a universal property, that of ZF-free algebra, whereas the class of all ordinals is characterized globally by the property of constituting ZF-free algebra with inflationary/monotonous successor function. In the same direction, the (small fibers) rank map, 
$\rho: V \to On$, is determined solely by the universal property of $V$, and the inclusion map, $i: On \to V $, is given by an adjunction condition.

The main aim of this  work is to obtain an "algebraic theory for surreals" along the lines of the Algebraic Set Theory: to establish 
abstract and general links between the class of all surreal numbers and a universe of "surreal sets" similar to (but expanding it) the (ZF-algebra) relations 
between the classes $On$ and $V$, giving the first steps towards  a certain  kind of  (alternative) "relative set theory" (see \cite{Fre16} for another presentation of this general notion).

In more details:

We want to perform a construction  (within the background class theory NBG) of a "class of all surreal sets", $V^*$, that satisfies, as far as possible, the following requirements:

\bu  $V^*$ is an expansion of the class of all sets $V$, via a map $j^* : V \to V^*$.

\bu $V^*$ is ranked in the class of surreal numbers $No$, in an analogous fashion that $V$  is ranked in the class of ordinal numbers $On$. I.e., expand,  through the injective map $j : On \to No$, the traditional set theoretic relationship $ V \overset{\rho}\to{\underset{i}\to{\rightleftarrows}} On$
to a new setting $V^* \overset{\rho^*}\to{\underset{i^*}\to{\rightleftarrows}} No$.

Noting that:\\
(i) the (injective) map $j : On \to No$, is a kind of "homomorphism", partially encoding the ordinal arithmetic;\\
(ii) the traditional set-theoretic constructions (in $V$) keep some relation with its (ordinal) complexity  (e.g., $x \in y \ \to\  \rho(x)< \rho(y)$, $\rho(\{x\}) = \rho(P(x)) = \rho(x)+1$, $\rho(\bigcup_{i \in I} x_i) = \bigcup_{i \in I} \rho(x_i)$);\\
then we wonder about the possibility of this new structured domain $V^*$ determines a category, by the encoding of arrows (and composition) as objects of $V^*$, in an analogous fashion to the way that the class $V$ of all sets determines a category, i.e. by the encoding of some notion of "function" as certain surreal set (i.e. an objects
of $V^*$); testing, in particular, the degree of compatibility of such constructions with the map $j^* : V \to V^*$ and examining if this new expanded domain could encode homomorphically the cardinal arithmetic.









We list below 3 instances of communications that we have founded in our bibliographic research on possible themes relating surreal numbers and set theory: we believe that they indicate  that we are pursuing a right track.

{\bf (I)} \underline{The Hypnagogic digraph and applications}\\
J. Hamkins have defined in \cite{Ham13}  the notion of "hypnagogic digraph", $(Hg, \rightharpoonup)$, an acyclic digraph graded on $(No,<)$, i.e., it is given a "rank" function $v : Hg \to No$ such that for each $x,y \in Hg$, if $x \rightharpoonup y$, then $v(x)< v(y)$.
 The hypnagogic digraph is a proper-class analogue the countable random $\mathbb{Q}$-graded digraph: it  is the Fra\"iss\'e limit of the class of all finite $No$-graded digraphs. It is simply the $On$-saturated $No$-graded class digraph,  making it set-homogeneous and universal for all class acyclic digraphs. 
Hamkins have applied this structure, and some relativized versions,  to prove interesting results concerning models of ZF set theory.


{\bf (II)} \underline{Surreal Numbers and Set Theory} \\
https://mathoverflow.net/questions/70934/surreal-numbers-and-set-theory\\
Asked July 21, 2011, by  Alex Lupsasca:\\
{\em \small I looked through MathOverflow's existing entries but couldn't find a satisfactory answer to the following question:\\
What is the relationship between $No$, Conway's class of surreal numbers, and $V$, the Von Neumann set-theoretical universe?\\
In particular, does $V$ contain all the surreal numbers? If so, then is there a characterization of the surreal numbers as sets in $V$? And does $No$ contain large cardinals?\\
I came across surreal numbers recently, but was surprised by the seeming lack of discussion of their relationship to traditional set theory. If they are a subclass of $V$, then I suppose that could explain why so few people are studying them.}


{\bf (III)} \underline{Surreal Numbers as Inductive Type?}\\
https://mathoverflow.net/questions/63375/surreal-numbers-as-inductive-type?rq=1\\
Asked in April 29, 2011, by Todd Trimble:\\
{\em \small  Prompted by James Propp's recent question about surreal numbers, I was wondering whether anyone had investigated the idea of describing surreal numbers (as ordered class) in terms of a universal property, roughly along the following lines.\\
In categorical interpretations of type theories, it is common to describe inductive or recursive types as so-called initial algebras of endofunctors. The most famous example is the type of natural numbers, which is universal or initial among all sets $X$ which come equipped with an element $x$ and an function $f : X \to X$. In other words, initial among sets $X$ which come equipped with functions $1+X \to X$ (the plus is coproduct); we say such sets are algebras of the endofunctor  $F$ defined by $F(X)=1+X$. Another example is the type of binary trees, which can be described as initial with respect to sets that come equipped with an element and a binary operation, or in other words the initial algebra for the endofunctor $F(X)=1+X^2$.\\
In their book Algebraic Set Theory, Joyal and Moerdijk gave a kind of algebraic description of the cumulative hierarchy $V$. Under some reasonable assumptions on a background category (whose objects may be informally regarded as classes, and equipped with a structure which allows a notion of "smallness"), they define a ZF-structure as an ordered object which has small sups (in particular, an empty sup with which to get started) and with a "successor" function. Then, against such a background, they define the cumulative hierarchy $V$ as the initial ZF-structure, and show that it satisfies the axioms of ZF set theory (the possible backgrounds allow intuitionistic logic). By tweaking the assumptions on the successor function, they are able to describe other set-theoretic structures; for example, the initial ZF-structure with a monotone successor gives $On$, the class of ordinals, relative to the background.\\
Now it is well-known that surreal numbers generalize ordinals, or rather that ordinals are special numbers where player $R$ has no options, or in different terms, where there are no numbers past the "Dedekind cut" which divides $L$ options from $R$ options. In any case, on account of the highly recursive nature of surreal numbers, it is extremely tempting to believe that they too can be described as a recursive type, or as an initial algebraic structure of some sort (in a background category along the lines given by Joyal-Moerdijk, presumably). But what would it be exactly?\\
I suppose that if I knocked my head against a wall for a while, I might be able to figure it out or at least make a strong guess, but maybe someone else has already worked through the details?}


\


{\bf Overview of the paper:}


\underline{Section 1:} \\
 This initial section establishes the notations and contains the preliminary  results needed for  the sequel of this thesis. It begins establishing our  set theoretic backgrounds -- that we will use freely in the text without further reference -- which is founded in NBG class theory, and contains mainly the definitions  and basic results on  some kinds of binary relations, in particular on  well-founded relations, and "cuts" as certain pairs of subsets of a class endowed with a binary irreflexive relation. After, in the second subsection, we introduce briefly a (categorically naive) version of ZF-algebras, a notion introduced in the 1990's in the setting of Algebraic Set Theory, in particular we introduce the concept of {\em standard ZF-algebra}, a concept also very useful to future developments (\cite{RM19}). The last subsection is dedicated to introduce the linearly ordered class of surreal numbers under many equivalent constructions and to present a characterization and some of its main structure, including its algebraic structure and its relations with the class (or ZF-algebra) of all ordinal numbers. 


\underline{Section 2:} \\
Motivated by properties of the linearly ordered class $(No,<)$, we introduce in this section the notion of  {\em Surreal Algebra (SUR-algebra)}: an structure ${\cal S} = (S,<, *, -, t)$, where $<$ is an acyclic relation on $S$, $*$ is a distinguished element of $S$, $-$ is an involution of $S$ and $t$ is a function that chooses an intermediary member between each small (Conway) cut in $(S,<)$,  satisfying some additional compatibility properties between them\footnote{Very recently, we came across with a study of surreal (sub)structures \cite{BH19}, that explores the theme under a different perspective.}. Every SUR-algebra turns out to be a proper class. We verify that $No$ provides naturally a SUR-algebra and present new relevant examples: the free surreal algebra ($SA$) and the free transitive SUR-algebra ($ST$).




\underline{Section 3:} \\
This section is dedicated to a generalization of this concept of SUR-algebra: we introduce the notion of partial SUR-algebra and consider two kinds on morphisms between them. This relaxation is needed to perform constructions as products, sub partial-SUR-algebra and certain  kinds of directed colimits. Some more examples are provided.

\underline{Section 4:} \\
As an application of the partial SUR-algebras theory previously worked out, we are able to prove in this section some universal properties satisfied by $SA$ and $ST$ (and generalizations), that justifies its  names of (relatively)  free SUR-algebras.

\underline{Section 5:} \\
In  this final section, we briefly comment on the sequel of this paper (\cite{RM19}) and present a  list of questions that have occurred to us during the elaboration of the work that we intend to address in the future.


\section{Preliminaries} \label{chap1}


 This section establishes the notations and contains the preliminary  results needed for  the sequel of this thesis. It begins establishing our  set theoretic backgrounds -- that we will use freely in the text without further reference -- which is founded in NBG class theory, and contains mainly the definitions  and basic results on  some kinds of binary relations. After, we introduce briefly a (categorically naive) version of ZF-algebras, a notion introduced in the 1990's in the setting of Algebraic Set Theory. Finally, we present the class of surreal numbers, and some of its main structure, under many equivalent constructions.



\subsection{Set theoretic backgrounds} \label{sec1.1}

This preliminary  section is devoted to establishing our  set theoretic backgrounds  which is  founded in NBG class theory\footnote{In some parts of the thesis, we will need some category-theoretic tools and reasoning, thus we will use an expansion of  NBG by axioms asserting the existence of  Grothendieck universes.}, and contains mainly the definitions   and basic results on  the binary relations that will appear in the sequel of this work as: (strict) partial order relations, acyclic  relations, extensional relations,   well founded relations, and "cuts" as certain pairs of subsets of a class endowed with a binary irreflexive relation.

\subsubsection{NBG} \label{NBG-sub}

In this work, we will adopt the (first-order, with equality) theory NBG as our background set theory, where the unique symbol in the language is the binary relation $\in$. We will use freely the results of NBG, in the sequel, we just recall below some notions and notations. We recall also the basic notions and results on some kinds of binary relations needed for the development of this thesis. Standard references of set/class theory are \cite{Jechbk} and \cite{Kunbk}.

\bct {On NBG:} \label{NBG-ct}
 
Recall that the primitive concept of NBG is the notion of {\em class}. A class  is {\em improper} when it is a member of some class, otherwise the class is {\em proper}. The notion of {\em set} in NBG is defined: a set is a improper class. 

We will use $V$ to denote the universal class -- whose members are all sets --; $On$ will stand for the class of all ordinal numbers and $Tr$ denote the class of all transitive sets. $On \sub Tr \sub V$ and all the tree are proper classes.

Given  classes ${\cal C}$ and ${\cal D}$,  then ${\cal C}$ is a subclass of ${\cal D}$ (notation: ${\cal C} \sub {\cal D}$), when all members of ${\cal C}$ are also members of ${\cal D}$. Classes that have the same members are equal. Every subclass of a set is a set.

$\empty$ stands  for the unique class without members. $\empty$ is a set.

Given a class ${\cal C}$, denote $P_s({\cal C})$ the class whose members are all the {\em subsets} of ${\cal C}$. If ${\cal C}$ is a set, then $P_s(\cal C)$  is a set too. There is no class that has as members all the {\em subclasses} of a proper class\footnote{This is a "metaclass" in NBG, i.e., an equivalence class of formulae in one free variable, modulo the $NBG$-theory: any such formula is not collectivizing.}.

Given  classes ${\cal C}$ and ${\cal D}$, and a function $f : {\cal C} \ra {\cal D}$, then the (direct) image $f[{\cal C}] = \{d \in {\cal D}: \exists c \in {\cal C}, d =f(c)\}$ is a subset of $D$, whenever $\cal C$ is a set.

Since NBG satisfies the axiom of global choice (i.e., there is a choice function on $V \setminus \{\empty\}$) and then every class (proper or improper) can be well-ordered, which implies nice cardinality results: as in ZFC, any set $X$  is equipotent to a unique cardinal number (= initial ordinal), called the its cardinality of $X$ (notation: $card(X)$); moreover, all the proper classes are equipotent -- we will denote $card(C) = \infty$ the cardinality of the proper class $C$ -- $\infty$ can be seen as a representation of the well-ordered the proper class $On$.
\qdr\ect



\bct {\bf Binary relations:} \label{bin-ct}

A relation $R$ is a class whose members are ordered pairs. The domain (respect., range) of  $R$  is the class of all first (respect., second) components of the ordered pair in the relation. The support of the relation $R$ (notation: $supp (R)$)  is the class obtained by the reunion of its domain and range. We will say that a binary relation is defined {\em on/over} its support class.

A relation $R$ is reflexive when $(x,x) \in R$ for each $x$ in the support of $R$; on the other hand, $R$ is irreflexive, when $(x,x) \notin R$ for each $x$ in the support of $R$. We will use $<, \prec, \triangleleft$ to denote general irreflexive relations; $\leq, \preceq, \sqsub$ will stand for reflexive relations. 
	
A pre-order is a reflexive and transitive relation. A partial order is a antisymmetric pre-order. A {\em strict} partial order is a irreflexive and transitive relation. There are well known processes of: obtain a strict partial order from a partial order and conversely.

Let $R$ be a binary relation and let $s, s' \in supp(R)$. Then $s$ and $s'$ are $R$-comparable when: $s=s'$ or $(s,s') \in R$ or $(s',s) \in R$. A relation $R$ is total or linear when every pair of members of its support are comparable.

Every pre-order relation $\preceq$ on a class ${\cal C}$  gives rise to an equivalence relation $\sim$ on the same support: for each $c, c' \in {\cal C}$, $c \sim c'$ iff  $c \preceq c'$ and $c' \preceq c$.

Let $n \in \N$, a $n$-cycle of the relation $R$ is a $n+1$-tuple $(x_0, x_1, \cdots, x_n)$ such that $x_n = x_0$ and, for each $i<n$, $(x_i, x_{i+i}) \in R$. A relation is {\em acyclic} when it does not have cycles. Every acyclic relation is irreflexive. A binary relation is a strict partial order iff it is a transitive and acyclic relation. Note that a binary relation is  acyclic and total iff it is a strict linear order. 
\qdr\ect

\bct {\bf Induced binary relations:} \label{binind-ct}

Given a binary relation $R$ on a class ${\cal C}$. For each $c \in {\cal C}$, denote $c^{R} := \{d \in C: (d,c) \in R\}$. 

Define a new binary relation on ${\cal C}$: for each $c, c' \in {\cal C}$, $c \sqsub_R c'$ iff holds $\all x ( (x,c) \in R \ra (x, c') \in R)$ or, equivalently, $c^R \sub {c'}^R$.
Clearly, $\sqsub_R$ is pre-order relation on ${\cal C}$.

Denote $\equiv_R$, the equivalence relation associated to the pre-order $\sqsub_R$ . We will say that the binary $R$ is extensional  when $\equiv_R$ is the identity relation on ${\cal C}$ or, equivalently, when $\sqsub_R$ is a partial order. The axiom of extensionality in NBG ensures that $(V, \in_{\rest V \times V})$ is a class endowed with an extensional relation and, since members of ordinal numbers are ordinal numbers\footnote{If $\alpha \in On$, then $\alpha^{\in} = \{\beta \in On: \beta \in \alpha\} = \{x \in V: x \in \alpha\} = \alpha$.), then $(On, \in_{\rest On \times On})$ is class endowed with an extensional relation. }
\qdr\ect


\subsubsection{Well founded relations}  \label{wellfounded-sub}

	In this subsection we recall basic properties and constructions  concerning general well-founded relations. Also, we introduce a special kind of well-founded relation  suitable for our purposes.

\bct {\bf On well-founded relations:} \label{wfr-ct}

Recall that a binary  relation $\prec$ on a class ${\cal C}$ is well-founded when:\\
(i) The subclass $x^{\prec} = \{y \in {\cal C}: y \prec x\}$ is a set.\\
(ii) For each $X \sub {\cal C}$ that is a non-empty subset, there is $u \in X$ that is a  $\prec$-minimal member of $X$ (i.e., $\all v \in {cal C}$, $v \prec u \Ra v \notin X$).\footnote{By the global axiom of choice  (for classes), this condition is equivalent of a apparently stronger one: \\
(ii)' For each $X \sub {\cal C}$ that is a non-empty {\em subclas}, there is $u \in X$ that is a  $\prec$-minimal member of $X$.}

Let $\prec$ be an well-founded relation on a class ${\cal C}$. Since for each $n \in \N$, the (non-empty) subset $\{x_0, \cdots, x_n\} \sub {\cal C}$ has a $\prec$-minimal member, then $\prec$  is an acyclic relation and, in particular, $\prec$ is irreflexive.

If ${\cal D} \sub {\cal C}$, then $({\cal D}, \prec_{\rest {\cal D} \times {\cal D}})$ is an well-founded class.

An well-founded relation that is a strict linear/total order is a well-order relation.

The axiom of regularity in NBG, guarantees that the binary relation $\in$ over the universal class $V$ is an well-founded relation. $(On, \in)$ is an well-ordered proper class.

Let $\prec$ be an well-founded relation on a class ${\cal C}$. Then it holds:\\
{\bf The induction principle:} Let $X \sub {\cal C}$ be a subclass. If, for each $c \in {\cal C}$, the inclusion $c^{\prec} \sub X$ entails $c \in {\cal C}$, then $X= {\cal C}$.\\
{\bf The recursion theorem:} Let $H$ be a (class) function such that $H(c, g)$ is defined for each  $c \in {\cal C}$ and $g$  a (set) function with domain $c^{\prec}$. Then there is a unique (class) function $F$ with domain ${\cal C}$ such that $F(c) = H(c, F_{\rest c^{\prec}})$, for each  $c \in {\cal C}$.
\qdr\ect

\bct {\bf Rooted well-founded relations:} \label{wfrrooted-ct}

{\bf Remark:} Let $({\cal C}, \prec)$ be a well-founded class; the subclass $root({\cal C})$ of its {\em roots} has as members its $\prec$-minimal members. Note that:\\
\bu If ${\cal C}$ is a non-empty class, then $root({\cal C})$ is a non-empty class.\\
\bu If $\sqsub$ denotes the pre-order on ${\cal C}$ associated to $\prec$ (i.e., $c\sqsub d$ iff $\forall x \in {\cal C} (x\prec c \Ra x \prec d)$, then: $root({\cal C}) = \{ a \in C : a \sqsub c $, for  all $c \in {\cal C}\} $.

{\bf Definition:} A well-founded class $({\cal C}, \prec)$ will be called {\em rooted}, when it has a unique root $\Phi$. 
If it is the case, then the structure $({\cal C}, \prec, \Phi)$ will be called a rooted  well-founded class.

If $\prec$ is an extensional  well-founded relation on a non-empty class ${\cal C}$, then $({\cal C}, \prec)$  is rooted: indeed, if $r, r' \in root({\cal C})$, then $r \sqsub r'$ and $r' \sqsub r$, thus $r=r'$. However, to emphasize  the distinguished element in a structure of rooted  well-founded class, we will employ the redundant expression "rooted extensional well-founded class".
	 
{\bf Examples and counter-examples:} 
	
	$(V, \in, \empty)$ is a rooted extensional well-founded class
	
 $(On, \in, \empty)$ is a rooted extensional well-ordered class.
	
	Every  well-ordered set $(X, \leq)$ gives rise to a rooted extensional well-ordered set $(X, <, \Phi)$, where $ \Phi =\bot$ is the least element of $X$ and the strict relation, $<$, associated to $\leq$, is an well-founded relation, since for each $x, y \in X$, $x^{<} \sub y^{<}$ iff $x \leq y$.
	
	$(\mathbb{N}, \leq, 0)$ is a well-ordered set, thus it gives rise to a rooted extensional well-ordered set. $(\mathbb{N} \setminus \{0\}, |, 1)$ is determines a rooted well-founded set that is not extensional. Note that $(\mathbb{N} \setminus \{0, 1\}, |)$ is an well-founded set that is not rooted  since its subset of minimal elements is the infinite set of all prime numbers.
\qdr\ect




\subsubsection{Cuts and densities} \label{Cut-sub}

Many useful variants of the concept of Dedekind cut were already been defined on the  setting  strict linear order on a given set (see for instance \cite{Allbk}). 
In this preliminary subsection  we present expansions of these notions in two different direction:  we consider binary relations that are only irreflexive (instead of being a strict linear order) and defined on  general classes instead of improper classes (= sets). We also generalize the notions of density a la Hausdorff to this new setting.

Through this subsection, ${\cal C}$ denote a class and $<$ stands for a irreflexive binary relation whose support is ${\cal C}$.

\bct {\bf Cuts} \label{cut-ct}

A {\bf Conway cut} in $({\cal C}, <)$ is a pair $(A,B)$ of {\em arbitrary subclasses}\footnote{$A$ and/or $B$ could be the empty set.} of ${\cal C}$ such that $\all a \in A, \all b \in B, a<b$ (notation $A<B$). Since $<$ is a irreflexive relation on ${\cal C}$, then $A \cap B = \empty$. A Conway cut $(A,B)$ will be called {\em small}, when $A$ and $B$ are {\em subsets} of ${\cal C}$.
We can define in theory NBG the class $C_s({\cal C}, <) := \{ (A,B) \in P_s({\cal C}) \times P_s({\cal C}) : A < B \}$, formed by all the small Conway cuts in $({\cal C}, <)$.

A {\bf Cuesta-Dutari cut} in  $({\cal C}, <)$ is a Conway cut $(A,B)$ such that $A \cup B = {\cal C}$. Note that $(\empty, {\cal C})$ and $({\cal C},\empty)$ are always Cuesta-Dutari cuts in $({\cal C}, <)$. On the other hand, if ${\cal C}$ is a proper class, then the class $CD_s({\cal C}, <)$ of all small Cuesta-Dutari cuts in $({\cal C}, <)$ is the empty class.

A {\bf Dedekind cut} in  $({\cal C}, <)$ is a Cuesta-Dutari cut $(A,B)$ such that $A$ and $B$ are non-empty subclasses. If ${\cal C}$ is a set, then $(A,B)$ is a Dedekind cut in  $({\cal C}, <)$ iff $(A,B)$ is a Conway Cut such that the set $\{A,B\}$ is a partition of ${\cal C}$.
\qdr\ect

\bct {\bf Densities}		\label{dense-ct}

Let $\alpha$ be an ordinal number. Then $({\cal C}, <)$ will be called an $\eta_\alpha$-class, when for each {\em small Conway cut} $(A,B)$ in  $({\cal C}, <)$, such that $card(A), card(B) < \aleph_\alpha$, there is some $t \in {\cal C}$ such that $\all a \in A, \all b \in B$, ($a<t$, $t<b$) (notation: $A<t<B$). 

Let $({\cal C}, <)$  be  an $\eta_\alpha$-class. Taking cuts $(\empty, \{c\})$ (respec. $(\{c\}, \empty)$), for all $c \in {\cal C}$, we can conclude that an $\eta_\alpha$-class $({\cal C}, <)$  does not have $<$-minimal (respec. $<$-maximal) elements. Taking cuts $(\empty, X)$ (or $(X, \empty)$), for all $X \sub {\cal C}$ such that $card(X) < \aleph_\alpha$, we see that an $\eta_\alpha$-class $({\cal C}, <)$  has $card({\cal C}) \geq \aleph_\alpha$.

An $\eta_0$-class $({\cal C}, <)$ is just  a "dense and  without extremes" class.

If $({\cal C}, <)$ is an $\eta_\alpha$-class and $\beta \in On$ is such that $\beta \leq \alpha$, then clearly $({\cal C}, <)$ is an $\eta_\beta$-class.

$({\cal C}, <)$ will be called an $\eta_\infty$-class, when it is an $\eta_\alpha$-class for all ordinal number $\alpha$: this means that for each {small Conway cut} $(A,B)$ in  $({\cal C}, <)$ there is some $t \in {\cal C}$ such that  $A<t<B$. Every $\eta_\infty$-class is a proper class. We will see  that the strictly linearly ordered proper class of all surreal numbers $(No,<)$ is  $\eta_\infty$. We will introduce in Section 2 the notion of SUR-algebra: every such structure is a $\eta_\infty$-class.
\qdr\ect	

{\underline{From now on, we will use the notion of Conway cut only in the {\em small} sense.}}


					

\subsection{On Surreal Numbers} \label{sec1.3}

This section is dedicated to present the class of surreal numbers -- a concept  introduced by J.H. Conway in the mid
1970's -- under many (equivalent) constructions within the background set class theory NBG, its order and algebraic structure and its relations with the class (or the ZF-algebra) of all ordinal numbers.



\subsubsection{Constructions}\label{subsec1.3.1}


\subsubsection{1.3.1.1 \ The Conway's construction formalized in NBG}


We begin with the Conway's construction following the appendix his book \cite{Conbk}, in wich he gave a more formal construction.

We start defining, recursively, the sets $G_\alpha$ in order to define class of games.
\begin{enumerate}[(i)]
\item $G_0=\{\pa\emptyset\emptyset\}$
\item $G_{\alpha}=\{\pa AB:A,B\subseteq \displaystyle\bigcup_{\beta<\alpha}G_\beta\}$
\end{enumerate}

The class $G$ of Conway games is given by $G=\bigcup_{\alpha<\on}G_\alpha$.

We define can put a preorder $\leqslant$ in $G$:

\begin{center}
$x\leqslant y$ iff not $x^L$ satisfies $x^L\geqslant y$ and not $y^R$ satisfies $x\geqslant y^R$.
\end{center}

The second step of the construction is the definition of the class of pre-numbers. We will again define the ordianl steps $P_\alpha$ recursively:

\begin{itemize}
\item $P_0=\{\pa\emptyset\emptyset\}$
\item $P_{\alpha}=\{\pa AB:A,B\in \displaystyle\bigcup_{\beta<\alpha}P_\beta\hbox{ e }B\not >A\}$
\end{itemize}

The class $P$ of the pre-numbers is given by $P=\displaystyle\bigcup_{\alpha<\on}P_\alpha$.

FInally, the class $No$ is defined as the quotient of the class of prenumbers by the equivalence relation induced by $\leqslant$. To avoid problems with equivalence classes that are proper classes, wue can make a Scott's Trick.

Following Conway's notation, we will denote a class $(X,Y)/\sim $ by $\{X | Y\}$ and given a surreal number $x$, we will denote $x=\{L_x|R_x\}$, where ${L_x|R_x}$ is a prenumber that represents $x$. We will also use the notation $x^L$ for an element of $L_x$ and $x^R$ for an element of $R_x$.

The birth function $b$ is defined as $b(x)=min\{\alpha:{\exists (L,R)\in P_\alpha} \  {x=\{L|R\}}\}$

We can also define, for any ordinal $\alpha$, the sets $O_\alpha$, $N_\alpha$ and $M_\alpha$ ("Old", "made" and "new"):

\noindent\bu $O_\alpha=\{x\in No:b(x)<\alpha\}$\\
\noindent\bu $N_\alpha=\{x\in No:b(x)=\alpha\}$\\
\noindent\bu $M_\alpha=\{x\in No:b(x)\leqslant\alpha\}$

To end this subsection we will now define, recursively, the operations $+,-,\cdot,\div$ in $P$:

\bu $x+y=\{x^L+y,x+y^L|x^R+y,x+y^R\}$;\\
\bu $xy=\{x^Ly+xy^L-x^Ly^L,x^Ry+xy^R-x^Ry^R|x^Ly+xy^R-x^Ly^R,x^Ry+xy^L-x^Ry^L\}$;\\
\bu $-x=\{-x^R|-x^L\}$;\\
\bu $0=\{\emptyset|\emptyset\}$;\\
\bu $1=\{0|\emptyset\}$.

\bpr With this operations, $No$ is a real-closed Field. In addition, every (set) field has an isomorphic copy inside $No$. If Global Choice is assumed, this is valid also for class fields.
\epr



\subsubsection{1.3.1.2 \  The Cuesta-Dutari cuts construction}

Given an order $(T,<)$, we can make a Cuesta "completion" of $T$, denoted by $\chi(T)$, wich is defined by

$$\chi(T)=(T\cup CD(T),<'),$$

with $<'$ defined as follows:

\begin{enumerate}[(i)]
  \item If $x,y\in T$, then the order is the same as in $T$;
  \item If $x=(L,R),y=(L',R')\in CD(T)$, then $x<'y$ if $L\subsetneq L'$;
  \item If $x\in T$ and $y=(L,R)\in CD(T)$, then $x<'y$ if $x\in L$ or $y<'x$ if $x\in R$.
\end{enumerate}

The idea of that construction is basically the iteration of the Cuesta-Dutari completion starting from the empty set until the we obtain a $\eta_{On}$ class.




By recursion we define the sets $T_\alpha$:
\begin{itemize}
  \item $T_0=\emptyset$;
  \item $T_{\alpha+1}=\chi(T_\alpha)$;
  \item $T_\gamma=\displaystyle\bigcup_{\beta<\gamma}T_\beta$.
\end{itemize}

And finally we have $${\no}=\displaystyle\bigcup_{\alpha\in \rm{\on}}T_\alpha$$

In that construction, the birth function $b$ is given by the map that assigns to each surreal number $x$, the ordinal $b(x)$ which corresponds to the set $T_b(x)$ that $x$ belongs.

Note that in this construction the sets "old", "made" and "new" can be presented in a simpler way:\\
\bu $O_\alpha=\displaystyle\bigcup_{\beta<\alpha} T_\alpha$\\
\bu $M_\alpha =\displaystyle \bigcup_{\beta\leqslant\alpha} T_\alpha$\\
\bu $N_\alpha=T_\alpha$


\subsubsection{1.3.1.3 \  The binary tree construction or the space of signs construction}





Consider the class $\Sigma=\{f:\alpha\to \{-,+\}:\alpha \in On\}$. We can define in this class an relation $<$ as follows:

\noindent \bu $f<g \iff f(\alpha)<g(\alpha)$, where $\alpha$ is the least ordinal such that $f$ and $g$ differs, with the convention $-<0<+$ ($f(\alpha)=0$ iff $f$ is not defined in $\alpha$).

With this relation, $\Sigma$ is a linearly ordered class isomorphic to $(No,\leqslant)$.

In this construction, the birth function is given by the map $b:\Sigma\to On$, $f\mapsto \dom f$. 






\subsubsection{The axiomatic approach} \label{axiomSur-sub}

It is an well-known fact that the notion of real numbers ordered field can be completely described (or axiomatized) as a certain structure --of complete ordered field -- and every pair of such kind of structure are isomorphic under a unique ordered field isomorphism (in fact, there is a unique ordered field morphism between each pair of complete  ordered fields and it is, automatically,  an isomorphism).  In this subsection, strongly based on section 3 of the chapter 4 in \cite{Allbk},  we present a completely analogous description for the ordered class (or  ordered field) of surreal numbers.

\bdf \label{suraxiom-df} A full class of surreal numbers is a structure ${\cal S} = (S,<,b)$ such that:
(i) $(S,<)$ is a strictly linearly ordered class;\\
(ii) $b: S \to On$ is a surjective function;\\
(iii) For each (small) Conway cut $(L,R)$ in $(S,<)$,  the class $I_S(L,R) = \{ x \in S : L<\{x\}<R\}$ is non-empty and its subclass $m_S(L,R) = \{ x \in I_S(L,R) : \all y \in S, b(y) < b(x) \ra y \notin I_S(L,R)\}$ is a singleton;\\
(iv) For each  Conway cut $(L,R)$ in $(S,<)$ and each ordinal number $\alpha$ such that $b(z)<\alpha$, $\all z \in L \cup R$, $b(\{L|R\}) \leqs \alpha$, where  $\{L|R\}$ its unique member of $m_S(L,R)$.
\qdr\edf

\bre \label{suraxiom-re} Let ${\cal S} = (S,<,b)$ be a full class of surreal numbers.\\
\bu  Condition (ii) above entails that $S$ is a {\em proper} class.\\
\bu Condition (iii) above guarantees that $(S,<)$ is a $\eta_\infty$-class.\\
\bu Since the order relation in $On$ is linear (is an well-order), according the notation in condition (iii), $m_S(L,R) = \{ x \in I_S(L,R) : \all y \in S,  y \in I_S(L,R) \ra b(x) \leq b(y)\}$.\\
\bu By condition (iv), $b(\{\empty|\empty\}) = 0$.
\qdr\ere

As mentioned in section 3 of chapter 4 in \cite{Allbk}, by results proven in  Conway's book \cite{Conbk}, the constructions of surreal numbers classes presented in our subsection 3.1 (by Conway cuts, by Cuesta-Dutari cuts and by the space of sign-expansions), endowed with natural "birthday" functions,  are all full classes of surreal numbers. It a natural question 

\bdf \label{surmono-df} Let ${\cal S} = (S,<,b)$ and ${\cal S}' = (S',<',b')$ be full classes of surreal numbers. 
A surreal (mono)morphism  $f : {\cal S} \to {\cal S}'$ is a function $f : S \ra S'$ such that:\\
(i) $\all x, y \in S, \ x< y \ra f(x) <' f(y)$;\\
(ii) $\all x \in S, \ b'(f(x))=b(x)$.
\qdr\edf

\bre \label{surmono-re} Let ${\cal S} = (S,<,b)$ and ${\cal S}' = (S',<',b')$ be full classes of surreal numbers. \\
\bu Since $<$ and $<'$ are linear order, a surreal morphism is always injective and condition (i) is equivalent to:\\
(i)'  $\all x, y \in S, \ x< y \lra f(x) <' f(y)$.\\
\bu Naturally, we can define a ("very-large") category whose objects are the full classes of surreal numbers and the arrows are surreal morphisms, with obvious composition and identities. Clearly, an isomorphism in such category is just a surjective morphism.
\qdr\ere

\bpr \label{surmono-pr} Let ${\cal S} = (S,<,b)$ and ${\cal S}' = (S',<',b')$ be full classes of surreal numbers. Then:\\
(i) There is a unique surreal (mono)morphism $f : {\cal S} \to {\cal S}'$ and it is an isomorphism. \\
(ii) For each ordinal number $\alpha$,  $b^{-1}([0,\alpha))$ is a set. Or, equivalently, $b$ is a locally small function.\\
(iii) The function $(L,R)\in C_s(S,<) \ \overset{t}\to\mapsto \ \{L|R\} \in S$ is surjective.
\qdr\epr

In particular, all the  constructions  of surreal numbers classes presented in our subsection 3.1, endowed with natural birthday functions,  are canonically isomorphic, through a unique isomorphism.

\

In the section 4 of chapter 4 in \cite{Allbk}, named "Subtraction in $No$", we can find the following result:

\bpr \label{subtraction-pr} Let ${\cal S} = (S,<,b)$ be a full class of surreal numbers. Then there is a unique function $- : S \to S$ such that:\\
(i) $b(-x) = b(x), \all x \in S$;\\
(ii) $-(-x) = x, , \all x \in S$;\\
(iii) $x<y \ra -y < -x, \all x, y \in S$;\\
(iv) $-\{L|R\} = \{-R|-L\}, \all (L,R) \in C_s(S,<)$.
\qdr\epr

\bre \label{subtraction-re} Let ${\cal S} = (S,<,b)$ be a full class of surreal numbers.\\
\bu In the presence of condition (ii), condition (iii) is equivalent to:\\
(iii)' $x<y \to -y < -x, \all x, y \in S$.\\
\bu By condition (iii), condition (iv) makes sense, since $L<R \Ra -R<-L$.\\
\bu By condition (iv), $-\{\empty|\empty\} = \{\empty|\empty\}$.
\qdr\ere

\

We finish this Subsection registering the following useful results
whose proofs can be found in \cite{Allbk}, pages 125, 126.

\bfa \label{Nocomp-fa} Let ${\cal S} = (S,<,b)$ be a full class of
surreal numbers. Let $A,A', B, B', \{x\}, \{x'\} \sub S$ be subsets
such that $A<B$ and $A'<B'$ and $x = \{A|B\}, x' = \{A'|B'\}$. Then:

(a) If $A$ and $A'$ are mutually cofinal and $B$ and $B'$ are mutually
coinitial, then $\{A|B\} = \{A'|B'\}$.

(b) Suppose that $(A,B)$ and $(A',B')$ are {\underline{timely
representations}} of $x$ and $x'$ respectively, i.e $b(z) < b(x), \all
z \in A \cup B$ and $b(z') < b(x'), \all z' \in A' \cup B'$. If $x =
x'$ then  $A$ and $A'$ are mutually cofinal and $B$ and $B'$ are
mutually coinitial.
\qdr\efa




\subsubsection{Ordinals in No}

The results presented in this Subsection can be found in the chapter 4 of \cite{Allbk}.

The ordinals can be embedded in a very natural way in the field $No$. The function that makes this work is recursively defined as follows:

\bdf
$j(\alpha) = \{j[\alpha]|\empty\}$, $\alpha \in On$.
\edf

The following result establishes a relation between the function $j$ and the birthday function:

\bpr
$b\circ j = id_On$
\epr

That map $j$ encodes completely  the ordinal order into the surreal order:

\bpr
$\alpha < \beta$ iff $j(\alpha) < j(\beta)$, $ \all \alpha, \beta \in On$.
\epr

We have also that $j(0) =0$, $j(1) = 1$. In fact, that embedding preserves also some algebraic structure. Although the sum and product of ordinals are not commutative, we have alternative definitions sum and product in $On$ closely related to the usual operations that are commutative:

\bdf
If $\alpha$ and $\beta$ are ordinals we can define the Hessemberg Sum of $\alpha$ and $\beta$
\edf

\bfa
The Hessemberg sum and products of ordinals are mapped by $j$ to the surreal sum and product.
\efa

In other words, the semi-ring $(On,\dot +,\dot \times,0,1)$ has an isomorphic copy in $No$ given by the image of $j$

\section{Introducing Surreal Algebras} \label{chap2}

Motivated by structure definable in the class $No$ of all surreal numbers, we introduce in this section the notion of surreal algebra (SUR-algebra) as a  (higher-order) structure $\cS = (S,<,*,-,t)$, satisfying some properties were, in particular, $<$ is an acyclic relation on $S$  where $t : C_s(S) \to S$ is a function that gives a coherent choice of witness of $\eta_\infty$ density of $(S,<)$. Every SUR-algebra turns out to be a proper class. Besides the verification of $No$  indeed support the SUR-algebra structure, we have defined two distinguished SUR-algebras $SA$ and $ST$, respectively the "free surreal algebra"' and "the free transitive surreal algebra" that will be useful in the sequel of this work. We also have introduced the notion of partial SUR-algebra (that can be a improper class)  and describe some examples and constructions in the corresponding categories.  We have provided, by categorical methods, some universal results that characterizes the SUR-algebras $SA$ and $ST$, and also some relative versions with base ("urelements") $SA(I)$, $ST(I')$ where $I$, $I'$ are partial SUR-algebra satisfying a few constraints. 

\subsection{Axiomatic definition} \label{sec2.1}

\bdf \label{Suralg-df} A surreal algebra (or SUR-algebra) is an structure $\cS = (S,<,*,-,t)$ where:\\
\bu  $<$ is a binary relation in $S$;\\
\bu $* \in S$ is a distinguished element;\\
\bu $- : S \to S$ is a unary operation;\\
\bu $t : C_s(S) \to S$ is a function, where $C_s(S) = \{ (A,B) \in P_s(S) \times P_s(S) : A<B\}$.

Satisfying the following properties:
\begin{enumerate}[{\bf (S1)}]
  \item $<$ is an acyclic relation.
  \item $\all x \in S$, $-(-x) = x$.
  \item $-*=*$.
  \item $\all a,b \in S$, $a<b$ iff $-b<-a$.
  \item $\all (A,B) \in C_s(S)$, $A<t(A,B)<B$.
  \item $\all (A,B) \in C_s(S)$, $-t(A,B)=t(-B,-A)$.
  \item $*=t(\emptyset,\emptyset)$.
\end{enumerate}
\qdr\edf

\bre \label{Suralg-re}\\
\bu  Let $(S, <)$ be the underlying  relational structure of a surreal algebra ${\cal S}$. Then $<$ is an  irreflexive relation, by condition (S1), and by (S5), $(S,<)$ is a $\eta_\infty$-relational structure. As a consequence $S$ is a {\em proper} class: \ref{dense-ct} in the Subsection \ref{Cut-sub}. The other axioms establish the possibility of choice of witness for the $\eta_\infty$ property satisfying additional coherent conditions.\\
\bu Note that (S3) follows from (S7) and (S6) : $-* = - t(\empty,\empty) = t(-\empty, -\empty) = t(\empty, \empty) = *$. \\
\bu Axiom (S7) establish that the SUR-algebra structure is "an extension by definitions" of a simpler (second-order) language: without a symbol for constant $*$.\\
\bu  In the presence of (S2), condition (S4) is equivalent to:\\
(S4)' $\all a,b \in S$, $a<b$ $\Ra$ $-b<-a$.\\
\bu By condition (S4), condition (S6) makes sense, since $A<B \Ra -B<-A$ (and if $A$, $B$ are sets, then  $-A$, $-B$ are sets).
\qdr\ere

A morphism of surreal algebras is a function that  preserves all the structure on the nose. More precisely:

\bdf \label{SURmor-df} Let ${\cal S} = (S, <, -, *, t)$ and ${\cal S}' = (S', <', -', *', t')$ be  SUR-algebras. A morphism of SUR-algebras $h : {\cal S} \rightarrow {\cal S}'$ is a function $h : S \to S'$ that satisfies the conditions below:
\begin{enumerate}
\item [{\bf (Sm1)}] $h(*) = *'$.
\item [{\bf (Sm2)}] $h(-a) = -'h(a)$, $\all a \in S$.
\item [{\bf (Sm3)}] $a<b \implies h(a) <' h(b)$, $\all a,b \in S$.
\item [{\bf (Sm4)}] $h(t(A,B))=t'(h[A],h[B])$, $\all (A,B)\in C_s(S)$.\footnote{Note that, by property  (Sm3), $(A,B)\in C_s(S) \implies (h[A],h[B])\in C_s(S')$, thus (Sm4)  makes sense.}
\end{enumerate}
\qdr\edf

\bdf \label{CATSUR-df} {\bf The category of SUR-algebras:}

We will denote by  $SUR-alg$  the ("very-large") category such that $Obj(SUR-alg)$ is the class of all  SUR-algebras 
and $Mor(SUR-alg)$ is the class of all partial SUR-algebras morphisms, endowed with obvious composition and identities.
\qdr\edf

\bre \label{CATSUR-re}

Of course, we have the same "size issue" in the categories $ZF-alg$, of all ZF-algebras, and in $SUR-alg$: each object can be (respect., is a) proper class, thus it cannot be represented in NBG background theory this "very large" category. The mathematical (pragmatical) treatment of this question , that we will adopt in the present work, is to assume a stronger background theory: NBG (or ZFC) and also three Grothendieck universes $U_0 \in U_1 \in U_2$. The members of $U_0$ represents "the sets"; the  members of $U_1$  represents "the classes"; the members of $U_2$ represents "the meta-classes". Thus a category $\cC$ is: (i) "small", whenever $\cC \in U_0$; (ii) "large", whenever $\cC \in U_1 \setminus U_0$; (iii) "very large", whenever $\cC \notin U_2 \setminus U_1$.  \qdr\ere




\subsection{Examples and constructions} \label{sec2.2}

\subsubsection{The surreal numbers as SUR-algebras} \label{SURNo-subsec}

The structure $(No, <, b)$ of  full surreal numbers class,  according the Definition \ref{suraxiom-df} in the Subsection \ref{axiomSur-sub},  induces a unique structure of SUR-algebra $(No,<,-,*,t)$, where:\\
\bu The function $t : C_s(No,<) \to No$ is such that $(A,B) \mapsto t(A,B) := \{A|B\}$;\\
\bu The distinguished element $* \in No$ is given by $* := \{\empty|\empty\}$;\\
\bu The function $- : No \to No$ is the unique function satisfying the conditions in Proposition \ref{subtraction-pr} and Remark \ref{subtraction-re}.

\

This SUR-algebra has two distinctive additional properties:\\
\bu \ $t$ is a surjective function;\\
\bu \ $<$ is a strict linear order (equivalently, since $<$ is acyclic, $<$ is a total relation).


\subsubsection{The free surreal algebra} \label{SURSA-subsec}

We will give now a new example of surreal algebra, denoted $SA$\footnote{The "A" in $SA$  is to put emphasis on {\bf a}cyclic.}, which is not a linear order and satisfies a nice universal property on the category of all surreal algebras (see Section 2.4). The construction, is based is based on a cumulative Conway's cuts hierarchy  over a family of binary  relations.  \footnote{Starting from the emptyset, and performing a cumulative construction based on Cuesta-Dutari completion of a linearly ordered set, we obtain $No$: see for instance \cite{Allbk}.}.

We can define recursively the family of {\bf sets} $SA_\alpha$ as follows:

Suppose that, for all $\beta<\alpha$,  we have constructed the sets $SA_\beta$ and $<_\beta$, binary relations on $SA_\beta$,  and denote $SA^{(\alpha)}=\bigcup_{\beta<\alpha}SA_\beta$ and $<^{(\alpha)}=\bigcup_{\beta<\alpha}<_\beta$ . Then, for $\alpha$ we define:

\begin{itemize}
\item $SA_\alpha=SA^{(\alpha)}\cup\{\pa AB : A,B\subseteq SA^{(\alpha)}\hbox{ and }A<^{(\alpha)}B\}$.\footnote{The expression $\pa AB$ is just an alternative notation for the ordered pair $(A,B)$, used for the reader's convenience.}
\item $<_\alpha=<^{(\alpha)}\cup\{(a,\pa AB),(\pa AB,b) :\pa AB\in SA_\alpha \setminus SA^{(\alpha)} \hbox{ and }a\in A,b\in B\}$.
\end{itemize}

\begin{itemize}
\item The class $SA$\footnote{Soon, we will see that $SA$ is a {\em proper} class.} is the union $SA := \bigcup_{\alpha\in On}SA_\alpha$.
\item $< := \bigcup_{\alpha\in On} <_\alpha$ is a binary relation on $SA$.
\end{itemize}

{\bf Fact:} Note that that: \\
(a) \ $SA^{(0)} = \empty$, $SA^{(1)} = SA_0 = \{ \pa \empty\empty \}$ and  $SA_1 = \{\pa \empty\empty, \pa \empty{\{\pa \empty\empty\}}, 
\pa {\{\pa \empty\empty\}}\empty \}$. By simplicity, we will denote $* := \pa \empty\empty = 0$, $-1 := \pa \empty{\{*\}}$ and $1 := \pa {\{*\}}\empty$ thus $SA_1 = \{0, -1, 1\}$.\\
(b) \ $<_0 = \empty$ and $<_1 = \{(-1, 0), (0, 1)\}$.\\
(c) \ $-1$ and $1$ are $<$-incomparable.\\
(d) \ $SA^{(\alpha)} \sub SA_\alpha$, $\alpha \in On$.\\ 
(e) \ $SA_\beta \sub SA_\alpha$, $\beta \leq \alpha \in On$.\\
(f) \ $SA^{(\beta)} \sub SA^{(\alpha)}$, $\beta \leq \alpha \in On$.\\
(g) \ $<^{(\alpha)} = <_{\alpha} \cap SA^{(\alpha)}\times SA^{(\alpha)}$,  $\alpha \in On$ (by the definition of $<_{\alpha}$).\\
(h) \ $<_\beta = <^{(\alpha)} \cap SA_{\beta}\times SA_{\beta}$, $\beta < \alpha \in On$.\\
(i) \ $<_\beta = <_{\alpha} \cap SA_{\beta}\times SA_{\beta}$, $\beta \leq \alpha \in On$ (by items (g) and (h) above).\\
(j) \ $<_\alpha = < \cap SA_{\alpha}\times SA_{\alpha}$,  $\alpha \in On$.\\
(k) \ $C_s(SA_{\alpha}, <_{\alpha}) = C_s(SA,<) \cap P_s(SA_{\alpha}) \times P_s(SA_{\alpha})$, $\alpha \in On$ (by item (j)).\\ 
(l) \ $C_s(SA^{(\alpha)}, <^{(\alpha)}) = C_s(SA,<) \cap P_s(SA^{(\alpha)}) \times P_s(SA^{(\alpha)})$, $\alpha \in On$.

\

We already have defined $<$ and $* (= \pa\emptyset\emptyset$) in $SA$, thus we must define  $t : C_s(SA) \ra SA$ and $- : SA \ra SA$ 
to complete the definition of the structure $SA$: this will be carry out by recursion on  well-founded relations on $SA$ and 
$C_s(SA)$\footnote{From now on, we will omit the binary relation on a class when it is clear from the setting.} that will be defined below.

\

For each $x\in SA$, we define its rank as  $r(x) := min\{\alpha \in On : x\in SA_\alpha\}$. Since for each $\beta < \alpha$, $SA_\beta \sub SA^{(\alpha)} \sub SA_\alpha$, it is clear that $r(x) = \alpha$ iff $x \in SA_\alpha \setminus SA^{(\alpha)}$. 

\

Following Conway (\cite{Conbk}, p.291), we can define for the $SA$ setting the notions of: "old members", "made members" and "new members". More precisely, for each ordinal $\alpha$:\\
\bu The set of {\bf old} members w.r.t. $\alpha$ is the subset of $SA$ of all members "born {\bf before} day $\alpha$".  $O(SA,\alpha) : = SA^({\alpha)}$;\\
\bu The set of {\bf made} members w.r.t. $\alpha$ is the subset of $SA$ of all members "born {\bf on or before} day $\alpha$".  $M(SA,\alpha) : = SA_{\alpha}$;\\
\bu The set of {\bf new} members w.r.t. $\alpha$ is the subset of $SA$ of all members "born {\bf on} day $\alpha$".  $N(SA,\alpha) : =  SA_\alpha \setminus SA^{(\alpha)}$.

 \

We will denote $x \prec y$ in $SA$ \ iff \ $r(x) < r(y)$ in $On$.

\

{\bf Claim 1:} $\prec$ is an well-founded relation in $SA$.

\bdm Let $y \in SA$ and let $\alpha := r(y)$. Given $x \in SA$, $r(x) < \alpha$ iff $x \in O(SA, \alpha) = SA^{(\alpha)}$. Therefore, the  subclass $\{x \in SA: x \prec y\}$ is a subset of $SA$. \\
Now let $X$ be a non-empty subset of $SA$. Then $r[X]$ is a non-empty subset of $On$ and let $\alpha := min r[X]$. Consider any $a \in r^{-1}[\{\alpha\}] \cap X$, then clearly $a$ is a  $\prec$-minimal member of $X$. 
\qdr

\

We have an induced "rank"  on the class  (of small $<$-Conway cuts) $C_s(SA) = \{(A,B) \in P_s(A) \times P_s(B) : A < B\}$ $R(A,B) := min\{\alpha \in On : A \cup B \sub SA^{(\alpha)}\}$.  We can also define a binary relation on the class $C_s(SA)$: \\
 $(A,B) \triangleleft (C,D)$ in $C_s(SA)$ \ iff \ $R(A,B) < R(C,D)$ in $On$.

\

{\bf Claim 2:} $\triangleleft$ is an well-founded relation in $C_s(SA)$.

\bdm Let $(C,D) \in C_s(SA)$ and let $\alpha := R(C,D)$. Given $(A,B) \in C_s(SA)$, $R(A,B) < \alpha$ iff $\exists \beta < \alpha, A \cup B  \sub O(SA, \beta) = SA^{(\beta)}$. Therefore, the  subclass $\{(A,B) \in C_s(SA): (A,B) \triangleleft (C,D) \}$ is a subset of $C_s(SA)$.\\ 
Now let $Y$ be a non-empty subset of $C_s(SA)$. Then $R[Y]$ is a non-empty subset of $On$ and let $\alpha := min R[Y]$. Consider any $(A,B) \in R^{-1}[\{\alpha\}] \cap Y$, then clearly $(A,B)$ is a  $\triangleleft$-minimal member of $Y$. 
\qdr

\

Let $H$ be the (class) function  $H(p, g)$ where, for each  $p = (C,D) \in C_s(SA)$ and $g$  a (set) function with domain $p^{\triangleleft} :=\{ (A,B) \in C_s(SA) : (A,B) \triangleleft p\}$, given by $H(p, g) := \pa CD$ (i.e., $H$ is just first coordinate projection). Then $H$ is a class function and we can define by $\triangleleft$-recursion a unique (class) function $t : C_s(SA) \ra SA$ by $t(p) = H(p, t_{\rest {p^{\triangleleft}}})$, i.e. $t(C,D) = \pa CD$. The range of $t$ is included in $SA$: since $A$ and $B$ are {\em subsets} of $SA$ such that $A <B$, there exists $\alpha \in On$ such that $A, B \sub SA^{(\alpha)}$; since $<$ is the reunion of the increasing compatible family of binary relations $\{ <_\beta: \beta \in On\}$, the have that $A <^{(\alpha)} B$, thus $\pa AB \in SA_\alpha \sub SA$.

\

{\bf Claim 3:} $\all \alpha \in On,  M(SA, \alpha) = C_s(O(SA,{\alpha}))$. Thus 
$N(SA, \alpha) = C_s(SA^{(\alpha)}) \setminus SA^{(\alpha)}$.

\bdm Since $M(SA, \alpha) = O(SA,{\alpha}) \cup C_s(O(SA,{\alpha}))$, we just have to prove that, $SA^{(\alpha)} \sub C_s(SA^{(\alpha)}$, for each $\alpha \in On$. Suppose that $SA^{(\beta)} \sub C_s(SA^{(\beta)}$ for each $\beta \in On$ such that $\beta < \alpha$.
By the assumption,  we have $SA^{(\alpha)} = \bigcup_{\beta <\alpha} SA_\beta  = \bigcup_{\beta <\alpha}  C_s(SA^{(\beta)}$. Since $SA^{(\beta)} \sub SA^{(\alpha)}$ and $<^{(\beta)} = <^{(\alpha)} \cap (SA^{(\beta)} \times SA^{(\beta)}) $\footnote{The non  trivial inclusion $<^{(\beta)} \supseteq <^{(\alpha)} \cap (SA^{(\beta)} \times SA^{(\beta)})$ holds since for every pair $(x,y)$ in the right side there are $\delta < \beta$ and $\gamma < \alpha$ (that we can assume $\gamma \geq \delta$) such that $(x,y) \in <_\gamma \cap SA_\delta \times SA_\delta = <_\delta \sub <_\beta$.}, we have  $C_s(SA^{(\beta)}) \sub C_s(SA^{(\alpha))}$, thus  $\bigcup_{\beta <\alpha}  C_s(SA^{(\beta)} \sub  C_s(SA^{(\alpha)}$. Summing up, we conclude that $SA^{(\alpha)} \sub C_s(SA^{(\alpha)}$ and the result follows by induction.
\qdr

\

{\bf Claim 4:} $C_s(SA)  = SA$ and  $t : C_s(SA) \ra SA$ is the identity map, thus, in particular, $t$ is a bijection.

\bdm
By items (k) and (l) in the Fact above, $C_s(SA,<) = \bigcup_{\alpha \in On} C_s(SA_\alpha, <_\alpha) = \bigcup_{\alpha \in On} C_s(SA^{(\alpha)}, <^{(\alpha)})$. By Claim 3 above, $SA_\alpha =  C_s(SA^{(\alpha)},<^{(\alpha)}), \all \alpha \in On$, thus $\bigcup_{\alpha \in On} C_s(SA^{(\alpha)}, <^{(\alpha)}) = \bigcup_{\alpha \in On} SA_\alpha = SA$. Summing up, we obtain $C_s(SA) = SA$. Then $t : C_s(SA) \ra SA$, given by $(A,B) \mapsto \pa AB$ is the identity map.
\qdr

\

For each $x \in SA$, denote $(L_x,R_x) \in C_s(SA)$ the unique representation of $x$: in fact, $x = (L_x,R_x)$.

\

{\bf Claim 5:} $r \circ t = R$.

\bdm The functional equation is equivalent to:\\
 $\all (A,B) \in C_s(SA)$, $\all \gamma \in On$, \ $\pa AB \in SA_\gamma$  \ iff \ $A \cup B \sub SA^{(\gamma)}$. \\
If  $(A,B) \in C_s(SA)$ and $A \cup B \sub SA^{(\gamma)}$ then, since $A < B$ we have $A <^{(\gamma)} B$, thus $\pa AB \in SA_\gamma$ by the recursive definition of $SA_\gamma$. On the other hand, if $(A,B) \in C_s(SA)$ and  $\pa AB \in SA_\gamma$, then  by Claim 3 above,  $(A,B) \in C_s(SA,<) \cap C_s(SA^{(\gamma)}), <^{(\gamma)}) = C_s(SA^{(\gamma)}), <^{(\gamma)})$, thus $A \cup B \sub SA^{(\gamma)}$. 
\qdr

\

{\bf Claim 6:} Let $(A,B) \in C_s(SA)$ and $\alpha \in On$, then: $\all a \in A, \all b \in B$, $r(a), r(b) < \alpha$  \ iff \ $r(t(A,B)) \leq \alpha$. In particular: $\all a \in A, \all b \in B$, $r(a), r(b) < r(t(A,B))$.

\bdm The equivalence is just a rewriting of the equivalence proved above:  \\
$A \cup B \sub SA^{(\alpha)}$ \ iff \ $\pa AB \in SA_\alpha$.
\qdr

\

{\bf Claim 7:} $\all x, y \in SA$, $x < y$ \ $\im$ \ $r(x) \neq r(y)$. In particular, the relation $<$ in $SA$ is irreflexive.

\bdm Suppose that there are $x,y \in SA$ such that $x<y$ and $r(x) = r(y) = \alpha \in On$. Thus $x, y \in SA_{\alpha} \setminus SA^{(\alpha)}$ and, since $x,y \in SA_\alpha$ and $(x,y) \in < $, we get $(x,y) \in  <_\alpha \setminus <^{(\alpha)}$.  Thus $(x,y) = (a,(L_y,R_y))$ for some $a \in L_y \sub SA^{(\alpha)}$ or $(x,y) = ((L_x,R_x),d)$ for some $d \in R_x \sub SA^{(\alpha)}$. In both cases we obtain $x=a \in SA^{(\alpha)}$ or $y=d \in SA^{(\alpha)}$, contradicting our hypothesis. \qdr

\

{\bf Claim 8:} Let $A, B \sub SA$ be subclasses such that $A<B$, then $r[A] \cap r[B] = \emptyset$.

\bdm
Suppose that $A<B$ and that there are $a \in A$ and $b \in B$ such that $r(a) = r(b) \in r[A] \cap r[B]$. Then  $a<b$ and $r(a) = r(b)$, contradicting the Claim 7 above.
\qdr

\

{\bf Claim 9:} Let $(A,B), (C,D) \in C_s(SA)$. Then  $\pa AB < \pa CD$ \ iff \ $\pa AB \in C$ (then $r(\pa AB)< r(\pa CD)$) \  or \ $\pa CD \in B$ (then $r(\pa CD)< r(\pa AB)$).

\bdm ($\Leftarrow$) If $\pa AB = c \in C$ and $r(\pa CD) = \alpha$, then $(c,\pa CD) \in <_\alpha \sub <$, thus $\pa AB < \pa CD$. The other case is analogous.\\
($\Rightarrow$) Suppose that $\pa AB < \pa CD$. By Claim 7 above we have $\alpha = r(\pa AB) \neq r(\pa CD) = \gamma$. If $\alpha < \gamma$ we have $SA_\alpha \sub SA^{(\gamma)}$ and $\pa CD \in SA_{\gamma} \setminus SA_{(\gamma)}$, thus $(\pa AB, \pa CD) \in <_\gamma \setminus <^{(\gamma)}$ and we have $\pa AB \in C$. If $\gamma < \alpha$ we conclude, by an analogous reasoning,  that $\pa CD \in B$.
\qdr

\

{\bf Claim 10:} Let $(A,B), (C,D) \in C_s(SA)$. Then  $A < \pa CD <B$ and $R(A,B) \leq R(C,D)$ \ iff \ $A \sub C$ and $B \sub D$. 

\bdm ($\Leftarrow$) Let $R(C,D) = \alpha$, then $\all c \in C, \all d \in D$, \ $(c,\pa CD), (\pa CD,d) \in <_\alpha \sub <$. If $A \sub C$ and $B \sub D$, then $A< \pa CD < B$ and $A \cup B \sub C \cup D \sub SA^{(\alpha)}$, i.e. $R(A,B) \leq \alpha = R(C,D)$.\\
($\Rightarrow$)  Let $A < \pa CD < B$ and suppose that there is $a \in A \setminus C$, then $\pa {L_a}{R_a} = a  < \pa CD$.  Since $\pa {L_a}{R_a} \notin C$, then by Claim 9 above, we have $\pa CD \in R_a$, thus:\\
 $R(C,D) =_{Claim 5} r(\pa CD) <_{Claim 6} = r(\pa {L_a}{R_a}) = r(a) < r(\pa AB) = R(A,B)$. Analogously, if $A < \pa CD < B$ and $B \setminus D \neq \empty$, we obtain $R(C,D) < R(A,B)$.
\qdr

\

{\bf Claim 11:} For each $(A,B) \in C_s(SA,<)$, $A< \pa AB < B$. In particular, $(SA,<)$ is a $\eta_infty$ {\em proper} class.

\bdm By Claim 7 above, $<$ is an irreflexive relation. By Claim 10 above, for each $(A,B) \in C_s(SA,<)$, $A< \pa AB < B$, thus $(SA,<)$ is a 
$\eta_infty$ class. It follows from  \ref{dense-ct} in the Subsection \ref{Cut-sub}, that $SA$ is proper class.
\qdr

\

{\bf Claim 12:} For each $(A,B) \in C_s(SA,<)$ and each $z \in SA$ such that $A< z < B$, then $r(t(A,B)) \leq r(z)$.

\bdm Suppose that the result is false and let $\alpha$ the least ordinal such that there are $(A,B) \in C_s(SA)$ and $z \in SA$ such that $A<z<B$, but $r(z) < r(t(A,B)) = R(A,B) = \alpha$: thus $\alpha >0$. By a simple analysis of the cases $\alpha$ ordinal limit and $\alpha$ successor, we can see that there are $A' \sub A, B' \sub B$ such that $R(A',B') = \alpha' < \alpha$ and $A'<z < B'$, contradicting the minimality of $\alpha$\footnote{Hint: in the case $\alpha = \gamma +1$, use Claim 7.}.
\qdr

Define, by recursion on the well-ordered proper class $(On,<)$, a function $s : On \ra SA$ by $s(\alpha) := \pa {s[\alpha]}{\empty}$, $\alpha \in On$.

\

{\bf Claim 13:} $r \circ s = id_{On}$. In particular, the function $r : SA \to On$ is surjective and $SA$ is a proper class.

\bdm
 We will prove the result by induction on the well-ordered proper class $(On, <)$.  Let $\alpha \in On$ and suppose that $r(s(\beta) = \beta$, for all ordinal $\beta < \alpha$. Then:\\
 $r(s(\alpha)) = r(\pa {s[\alpha]}{\empty}) = r(t(s[\alpha],{\empty})) =_{Claim \ 5} R(s[\alpha], \empty) = min\{\gamma \in On : s[\alpha] \cup \empty \sub SA^{(\gamma)}\}$.\\
 By the induction hypothesis, we have:\\
(IH) \ $s(\beta) \in SA_{\beta} \setminus SA^{(\beta)}$, for all ordinal $\beta < \alpha$.\\
Since $s(\beta) \in SA_{\beta}$, we have $s(\beta) \in SA^{(\alpha)}$, $\all \beta < \alpha$.  If $s[\alpha] \cup \empty \sub SA^{(\gamma)}$ for some  $\gamma < \alpha$,  then $s(\gamma) \in SA^{(\gamma)}$, in contradiction with (IH). Summing up, we conclude that $r(s(alpha) = \alpha$, and the result follows by induction.
\qdr

\

{\bf Claim 14:} There is a unique function $ - : SA \to SA$,  such that:\\
(i) $\all x \in SA, r(-x) = r(x)$;\\
(ii) $\all x \in SA, -(-x) = x$;\\ 
(iii) $\all x, y \in SA$,  $x<y$ iff $-y <-x$;\\
(iv) $\all (A,B) \in C_s(SA), -t(A,B) = t({-B},{-A})$.

\bdm 
Let $z \in SA$ and suppose that a function $-$ is defined for all $x,y \in SA$, such that $x,y \prec z$, satisfying the conditions (i)--(iv) adequately restricted to the subset $SA^{(\alpha)}$, for $\alpha := r(z) \in On$. Then $z = \pa {L_z}{R_z} = t(L_z,R_z) \in SA_\alpha \setminus SA^{(\alpha)}$ for a unique $(L_z, R_z) \in C_s(SA)$ (Claim 4) and $\alpha$ is the least $\gamma \in On$ such that $L_z  \cup R_z \sub SA^{(\gamma)}$ (Claim 5), thus $\all x \in R_z \cup L_z, r(x) < r(z)$. Then $-x$ is defined $\all x \in L_z \cup R_z$, satisfying the conditions (i)--(iv) restricted to the subset $SA^{(\alpha)}$. Since $x<y$, $\all x \in L_z \all y \in R_z$, it holds, by condition (iii), $-y<-x$ then $-R_z < -L_z$ and since $-R_z, -L_z$ are the images of a function on sets, $(-R_z,-L_z) \in C_s(SA)$. Moreover, by condition (i), $\alpha$ is the least $\gamma \in On$ such that $-R_z \cup -L_z \sub SA^{(\gamma)}$, i.e. $t(-R_z,-L_z) = \pa {-R_z}{-L_z} \in SA_\alpha \setminus SA^{(\alpha)}$. Define $-z := t(-R_z,-L_z)$.\\
Now we will prove that the conditions (i)--(iv) still holds for all members in $SA_\alpha \supsetneq SA^{(\alpha)}$.\\
(i) Let $x \in SA_\alpha$. If $x \in SA_{(\alpha)}$, this condition holds by hypothesis. If $x \in  SA_\alpha \setminus SA^{(\alpha)}$, then by the recursive definition above, $-x = \pa {-R_x}{-L_x} \in SA_\alpha \setminus SA^{(\alpha)}$, thus $r(-x) = \alpha = r(x)$. Thus (i) holds in $SA_\alpha$.\\
(ii) Let $x \in SA_\alpha \setminus SA^{(\alpha)}$, then $-x, --x \in SA_\alpha \setminus SA^{(\alpha)}$ (by the validity of condition (i) on $SA_\alpha$ established above). $-(-x) = -(- t(L_x,R_x)) = -  t(-R_x,-L_x) = t(-(-L_x), -(-R_x)) = t(L_x, R_x) = x$, since by hypothesis the conditions (iii) and (ii) holds for members of $SA^{(\alpha)}$.  Thus (ii) holds in $SA_\alpha$.\\
(iii) We suppose that $\all x, y \in SA^{(\alpha)}$,  $x<y$ iff $-y <-x$. Let $x,y \in SA_\alpha$ such that $x<y$. If both $x, y \in SA^{(\alpha)}$ then, by hypothesis $-y<-x$. Otherwise, by Claim 7, there is exactly one between $x$, $y$ that is a member of  $SA_\alpha \setminus SA^{(\alpha)}$. By Claim 9:  if $r(y)< r(x) =\alpha$ then $y \in R_x$; if $r(x) < r(y) = \alpha$ then $x \in L_y$. Thus:  if $r(y)< r(x) =\alpha$ then $-y \in -R_x = L_{-x}$, thus $-y<-x$; if $r(x) < r(y) = \alpha$ then $-x \in -L_y = R_{-y}$, thus $-y<-x$. Then we ahve proved that  $\all x, y \in SA_{\alpha}$,  $x<y$ $\Ra$ $-y <-x$. Since the conditions (i) and (ii) have already be established on $SA_\alpha$, we also have  $\all x, y \in SA_{\alpha}$,  $-y<-x$ $\Ra$ $x= -(-x) <-(-y) = y$.\\
(iv) Suppose that  $t(A,B) = \pa AB = \pa {-B}{-A} = t(-B,-A)$ holds for all $(A,B) \in \bigcup_{\beta < \alpha} C_s(SA) \cap P_s(SA^{(\beta)} \times P_s(SA^{(\beta)} = \bigcup_{\beta < \alpha} C_s(SA^{(\beta)}) =  \bigcup_{\beta < \alpha} SA_{\beta} = SA^{(\alpha)}$. We must prove that the condition still holds for all $(C,D) \in C_s(SA) \cap P_s(SA^{(\alpha)} \times P_s(SA^{(\alpha)} = C_s(SA^{(\alpha)} = SA_{\alpha}$. Let $z= (C,D) = \pa CD = t(C,D) \in  SA_\alpha \setminus SA^{(\alpha)}$. Then just by the recursive definition of $-z$, we have $-z = - t(A,B) = t(-B,-A)$, as we wish.
\qdr

\

Finally, we will prove that $SA$ satisfies all the 7 axioms of SUR-algebra:

\begin{enumerate}

\item[(S7)] $*=t(\emptyset,\emptyset)$. \\
This holds by our definition of $*$.

\item[(S5)] $\all (A,B) \in C_s(SA)$, $A<t(A,B)<B$.\\
This holds by the Claim 10  above.

\item[(S6)] $\all (A,B) \in C_s(SA)$, $-t(A,B)=t(-B,-A)$.\\
This holds by Claim 14.(iv).

\item[(S3)] $-*=*$.\\
Since $-* = - t(\empty,\empty) = t(-\empty, -\empty) = t(\empty, \empty) = *$.

\item[(S2)] $\all x \in SA$, $-(-x) = x$.\\
This holds by Claim 14.(ii).

\item[(S4)] $\all a,b \in SA$, $a<b$ iff $-b<-a$.\\
This holds by Claim 14.(iii).

\item[(S1)] $<$ is an acyclic relation.\\
Suppose that $<$ is not acyclic and take $x_0<...<x_n<x_0$ a cycle in $(SA,<)$ of minimum length $n \in \N$. Since $<$ is an irreflexive relation (see the Claim 7), $n>0$. Let $\alpha = max\{r(x_i) : i \leq n\}$ and let $j$ be the least $i \leq n$ such that $r(x_j) = \alpha$.\\
\underline{If $j = 0$}: Since $x_0<x_1$ and $x_n < x_0$, then by Claim 7, $r(x_1), r(x_n)< r(x_0)$. Writing $x_0=\pa {L_{x_0}}{R_{x_0}}$ (since, by Claim 4, $SA = C_s(SA)$), we obtain from Claim 9 that $x_n \in L_{x_0}$ and $x_1 \in R_{x_0}$. As ${L_{x_0}}<{R_{x_0}}$, we have $x_n<x_1$ and then  $x_1<...<x_n<x_1$ is a cycle of length $n-1 <n$, a contradiction.\\
\underline{If $j > 0$}: Then define $j^{-} := j-1$ and $j^+ := j+1$ (respec. $j^+ = 0$), if $j < n$ (respec. $j=n$).\\
Then by Claim 7, $r(x_{j^-}), r(x_{j^+})< r(x_j)$ and by Claim 9: $x_{j^{-}} \in L(x_j)$ and  $x_{j^+} \in R_{x_j}$. As ${L_{x_j}}<{R_{x_j}}$, we have $x_{j^-} < x_{j^+}$ and  then we can take a sub-cycle of the original one omitting $x_j$: this new cycle has of length $n-1<n$, a contradiction.

\end{enumerate}






					
\subsubsection{The free transitive surreal algebra} \label{SURST-subsec}

We will give now a new example of surreal algebra, denoted $ST$\footnote{The "T" in $ST$  is to put emphasis on {\bf t}ransitive.}, which is a strict partial order\footnote{Recall that a binary relation is that is a strict partial order iff it is a transitive and  acyclic relation.} that is not linear and satisfies a nice universal property on the category of all {\bf transitive} surreal algebras (see Section 2.4). The construction is similar to the construction  of $SA$ in the previous subsection: it is based on a cumulative Conway's cuts hierarchy  over a family of binary (transitive) relations.

We can define recursively the family of {\bf sets} $ST_\alpha$ as follows:

Suppose that, for all $\beta<\alpha$,  we have constructed the sets $ST_\beta$ and $<_\beta$, binary relations on $ST_\beta$,  and denote $ST^{(\alpha)}=\bigcup_{\beta<\alpha}ST_\beta$ and $<^{(\alpha)}=\bigcup_{\beta<\alpha}<_\beta$ . Then, for $\alpha$ we define:

\begin{itemize}
\item $ST_\alpha=ST^{(\alpha)}\cup\{\pa AB : A,B\subseteq ST^{(\alpha)}\hbox{ and }A<^{(\alpha)}B\}$.
\item $<_\alpha=$ the transitive closure of the relation $<'_\alpha$, where $<'_\alpha := (<^{(\alpha)}\cup\{(a,\pa AB),(\pa AB,b):\pa AB\in ST_\alpha \setminus ST^{(\alpha)}\hbox{ and }a\in A,b\in B\})$.
\end{itemize}

\begin{itemize}
\item The (proper) class $ST$ is the union $ST := \bigcup_{\alpha\in On} ST_\alpha$.
\item $< := \bigcup_{\alpha\in On} <_\alpha$ is a binary (transitive) relation on $ST$.
\end{itemize}

\

The following result is straightforward an completely analogous to the corresponding items in the Fact in the previous Subsection on $SA$:

{\bf Fact 1:} Note that that: \\
(a) \  $ST^{(0)} = \empty$, $ST^{(1)} = ST_0 = \{ \pa \empty\empty \}$. By simplicity, we will denote $0 := \pa \empty\empty$, $1:= \pa \empty{\{0\}}$, $-1:= \pa {\{0\}}\empty$. Thus:  $ST_{0} = \{0\}$, $SA_1 = \{0,1,-1\}$.\\
(b) \ $<_0 = \empty$, $<_1 = \{(-1, 0), (0, 1), (-1,1)\}$.\\
(c) \ $-1<0<1$, $- 1<\pa{\{-1\}}{\{1\}}<1$, but $0, \pa{\{-1\}}{\{1\}}$ are $<$-incomparable.\\
(d) \ $ST^{(\alpha)} \sub ST_\alpha$, $\alpha \in On$.\\ 
(e) \ $ST_\beta \sub ST_\alpha$, $\beta \leq \alpha \in On$.\\
(f) \ $ST^{(\beta)} \sub ST^{(\alpha)}$, $\beta \leq \alpha \in On$.
\qdr

\

Analogously to in the $SA$ case,  we can define rank functions $r : ST \ra On$\footnote{For each $x \in ST$, $r(x) = \alpha \in On$ iff $x \in ST_\alpha \setminus ST^{(\alpha)}$.} and $R : C_s(ST) \ra On$ that induces well-founded relations on $ST$ and on $C_s(ST)$. 

\

The results below are almost all (the exception are the items (m), (n), (o)) analogous to corresponding items in the Fact in the previous Subsection on $SA$. However, the techniques needed in the proofs are different than in $SA$ case and deserve a careful presentation.

{\bf Fact 2:}  \\
(g) \ $<^{(\alpha)} = <_{\alpha} \cap SA^{(\alpha)}\times SA^{(\alpha)}$,  $\alpha \in On$.\\
(h) \ $<_\beta = <^{(\alpha)} \cap SA_{\beta}\times SA_{\beta}$, $\beta < \alpha \in On$.\\
(i) \ $<_\beta = <_{\alpha} \cap SA_{\beta}\times SA_{\beta}$, $\beta \leq \alpha \in On$.\\
(j) \ $<_\alpha = < \cap SA_{\alpha}\times SA_{\alpha}$,  $\alpha \in On$.\\
(k) \ $C_s(ST_{\alpha}, <_{\alpha}) = C_s(ST,<) \cap P_s(ST_{\alpha}) \times P_s(ST_{\alpha})$, $\alpha \in On$.\\ 
(l) \ $C_s(ST^{(\alpha)}, <^{(\alpha)}) = C_s(ST,<) \cap P_s(ST^{(\alpha)}) \times P_s(ST^{(\alpha)})$, $\alpha \in On$.\\
(m) $\all \alpha \in On$, $<_\alpha$ is a transitive and a acyclic relation on $ST_\alpha$. \\
(n)  $<$ is a transitive and  acyclic relation (or, equivalently, it is a strict partial order) on $ST$.\\  
(o) Let $x, y \in ST$ and denote $\alpha := max\{r(x), r(y)\}$. Then are equivalent:\\
\ \bu \ $x < y$.\\ 
\ \bu \ Exists $n \in \mathbb{N}$, exists $\{z_0, \cdots, z_{n+1}\} \sub ST_\alpha$ such that: $x = z_0, y = z_{n+1}$;  $z_j \in L_{z_{j+1}}$ or $z_{j+1} \in R_{z_j}$, for all $j \leq n$; $\{z_1, \cdots, z_n\} \sub ST^{(\alpha)}$.

\bdm Item (i) follows from items (g) and (h). Items (k) and (l) follows from item (j). Items (n) and (o) are  direct consequences of item (m), since $< = \bigcup_{\alpha \in On} <_\alpha$. 

(g) Clearly $<^{(\alpha)} \ \sub \ <_{\alpha} \cap SA^{(\alpha)}\times SA^{(\alpha)}$. To show the converse inclusion let $x, y \in SA^{(\alpha)}$ be such that $x<_\alpha y$ and let $x=x_0 <'_\alpha...<'_\alpha x_n = y$ be a sequence in $(ST_\alpha,<'_{\alpha})$ with the number $k = card(\{i\leqs n: r(x_i) = \alpha\}$ being minimum. We will show that $k =0$, thus the sequence is just $x=x_0 <^{(\alpha)}...<^{(\alpha)} x_n = y$ and then $x <^{(\alpha)} y$ because $<^{(\alpha)}$ is a transitive relation (since $<_\beta$, $\beta \in On$ is a transitive relation, by construction). Suppose, by absurd, that $k>0$ and let $j$ be the least $i \leq n$ such that $r(x_j) = \alpha$. By our hypothesis on $x,y$ we have $0<j<n$. Since $x_{j-1} <'_{\alpha} x_j <'_{\alpha} x_{j+1}$, we have $r(x_{j-1}), r(x_{j+1})< r(x_j)= \alpha$ and $ x_{j{-}1} \in L(x_j)$,  $x_{j+1} \in R_{x_j}$. As ${L_{x_j}}<^{(\alpha)}{R_{x_j}}$, we have $x_{j-1} <^{(\alpha)} x_{j+1}$ and  then we can take a sub-cycle of the original one omitting $x_j$: this new cycle has $k-1<k$ members with rank $\alpha$, a contradiction.

(h) We only prove the non-trivial inclusion. Let $x, y \in SA_{\beta}$ be such that $x<^{(\alpha)} y$. Since $<^{(\alpha)} = \bigcup_{\gamma < \alpha} <_\gamma$, let $\beta'$ be the least $\gamma < \alpha$ such that $x<_{\beta'}y$. We will prove that $\beta' \leq \beta$, thus we obtain $x<_\beta y$, as we wish.  Suppose, by absurd, that $\beta'>\beta$. Then $(x,y) \in <_{\beta'} \cap SA^{(\beta')} \times SA^{(\beta')}$, and by the item (g) proved above $(x,y) \in <^{(\beta')}$. Thus there is some $\gamma < \beta'$ such that $x<_\gamma y$, contradicting the minimality of $\beta'$.

(j) Let $x, y \in SA_{\alpha}$ be such that $x<y$. Since $< = \bigcup_{\gamma \in On} <_\gamma$, let $\alpha'$ be the least $\gamma \in On$ such that $x<_{\alpha'}y$. We will prove that $\alpha' \leq \alpha$, thus we obtain $x<_\alpha y$, as we wish.  Suppose, by absurd, that $\alpha'>\alpha$. Then $(x,y) \in <_{\alpha'} \cap SA^{(\alpha')} \times SA^{(\alpha')}$, and by the item (g) proved above $(x,y) \in <^{(\alpha')}$. Thus there is some $\gamma < \alpha'$ such that $x<_\gamma y$, contradicting the minimality of $\alpha'$.

(m) 
 By definition of $<_\gamma$, $<_\gamma$ is a transitive relation,  $\all \gamma \in On$.\\
Suppose that the statement is false and let $\alpha \in On$ be the least ordinal such that $(ST_\alpha, <_\alpha)$ has some cycle. Then  $\all \beta < \alpha$, $<_\beta$ is an acyclic relation but $<_\alpha$ has some cycle (or, equivalently, $<'_\alpha$ has some cycle). Let  $x_0 <'_\alpha...<'_\alpha x_n <'_\alpha x_0$ be a cycle in $(ST_\alpha,<'_{\alpha})$ with the number $k = card(\{i\leqs n: r(x_i) = \alpha\}$ being minimum. Note that $k>0$, otherwise $x_0, ..., x_n \in SA^{(\alpha)}$ and the cycle is $x_0 <^{(\alpha)}...<^{(\alpha)} x_n <^{(\alpha)} x_0$, thus there  is a $\beta <\alpha$ and a cycle $x_0 <_{\beta}...<_{\beta} x_n <_{\beta} x_0$ in $(ST_\beta, <_\beta)$, contradicting our hypothesis. \\
Let $j$ be the least $i \leq n$ such that $r(x_j) = \alpha$.\\
\underline{If $j = 0$}: Since $x_0 <'_\alpha x_1$ and $x_n <'_\alpha x_0$, then  $r(x_1), r(x_n)< r(x_0)=\alpha$. Writing $x_0=\pa {L_{x_0}}{R_{x_0}}$,  we have that $x_n \in L_{x_0}$ and $x_1 \in R_{x_0}$. As ${L_{x_0}}<^{(\alpha)}{R_{x_0}}$, we have $x_n<^{(\alpha)} x_1$,  and then  $x_1<'_\alpha...<'_\alpha x_n<'_\alpha x_1$ is a cycle in $(ST_\alpha, <'_\alpha)$ with $k-2<k$ members with rank $\alpha$, a contradiction.\\
\underline{If $j > 0$}: Then define $j^{-} := j-1$ and $j^+ := j+1$ (respect. $j^+ = 0$), if $j < n$ (respect. $j=n$).\\
Then $r(x_{j^-}), r(x_{j^+})< r(x_j)= \alpha$ and:  $ x_{j^{-}} \in L(x_j)$,  $x_{j^+} \in R_{x_j}$. As ${L_{x_j}}<^{(\alpha)}{R_{x_j}}$, we have $x_{j^-} <^{(\alpha)} x_{j^+}$ and  then we can take a sub-cycle of the original one omitting $x_j$: this new cycle has $k-1<k$ members with rank $\alpha$, a contradiction.
\qdr

\

Since the harder part was already done, we just sketch the construction of the SUR-algebra  structure $(ST,<,-,*,t)$:\\
\bu \ As in the $SA$ case, from the well founded relation on $C_s(ST)$ we can define recursively a function with range $ST$,  $t : C_s(ST) \ra ST$ by $t(A,B) = \pa AB$. We can prove, by induction, that $ST_\alpha = C_s(ST^{(\alpha)}$, $\alpha \in On$. Thus $t$ is a bijection (is the identity function). Moreover, if  $(A,B) \in C_s(ST)$, then $A< t(A,B) < B$.\\
\bu \ We define $* := 0 = t(\empty, \empty)$. \\
\bu \ As in the $SA$ case, we can define (recursively) the function $ - : ST \to ST$ by $- \pa AB := \pa {-B}{-A}$.

\

The verification of the satisfaction of the SUR-algebra axioms (S2)--(S7) are analogous as in the $SA$ case. The satisfaction of (S1) was proved in item (m) of Fact 2 above.




\subsubsection{The cut surreal algebra}

In this Subsection we present a generalization of the SA, ST constructions.
Given a surreal algebra $S$, we can define a new surreal algebra whose domain is $C_s(S)$ with the following relations and operations:

\bdf Let $(S,<,-,*,t)$ be a surreal algebra. Consider the following structure in $C_s(S)$\\
\bu $*'=(\emptyset,\emptyset)$\\
\bu $-'(A,B)=(-B,-A)$\\
\bu $(A,B)<'(C,D) \iff t(A,B)<t(C,D)$\\
\bu $t'(A,B)=(t[A],t[B])$\\
\edf

\bpr With this operations $(C_s(S),<',-',*',t')$ is a surreal algebra.
\epr

\bdm
\begin{enumerate}[{\bf (S1)}]
    \item $<'$ is acyclic because any cycle $(A_0,B_0)<'...<'(A_n,B_n)$ induces a cycle $t(A_0,B_0)<...<t(A_n,B_n)$ in $S$, which is acyclic.
    \item $-'-'(A,B))=-'(-B,-A)=(--A,--B)=(A,B)$.
    \item $-'*'=-'(\emptyset,\emptyset)=(-\emptyset,-\emptyset)=(\emptyset,\emptyset)$.
    \item $(A,B)<'(C,D)$ iff $t(A,B)<t(C,D)$ iff $-t(C,D)<-t(A,B)$ iff $t(-D,-C)<t(-B,-A)$ iff $(-D,-C)<'(-B,-A)$ iff $-'(C,D)<'-'(A,B)$.
    \item Let $(A,B)\in C_s(C_s(S))$. Then $A<'B$ and thus $t[A]<t[B]$. Since $S$ satisfies {\bf (S5)}, $t[A]<t(t[A],t[B])<t[B]$. By the definition of $<'$, $A<'(t[A],t[B])<'B$ and then $A<'t'(A,B)<'B$.
    \item $-'t'(A,B)=-'(t[A],t[B])=(-t[B],-t[A])=(t[-'B],t[-'A])=t'(-'B,-'A)$
    \item $t'(\emptyset,\emptyset)=(t[\emptyset],t[\emptyset])=(\emptyset,\emptyset)=*'$
\end{enumerate}
\edm\qdr

Some properties of $X$ are transferred to $C_s(X)$ as we can see in the above proposition:
\

\bpr \label{casesCut-pr}
\begin{enumerate}[(a)]
    \item If $X$ is transitive then $C_s(X)$ is transitive.
    \item If $X$ is linear then $C_s(X)$ is pre-linear, i.e., denote $\sim_t$ the equivalence relation on $C_s(S)$ given by $(A,B) \sim_t (C,D)$ iff $t(A,B) = t(C,D)$. Then it holds exactly one between of the alternatives: $(A,B) <' (C,D)$; $(A,B) \sim_t (C,D)$; $(C,D) <' (A,B)$.
\end{enumerate}
\epr

\bdm 
\begin{enumerate}[(a)]
\item Suppose that we have $(A_1,B_1), (A_2,B_2), (A_3,B_3)\in C_s(X)$ satisfying $(A_1,B_1)<' (A_2,B_2)<' (A_3,B_3)$. Then, by definition, $t(A_1,B_1)< t(A_2,B_2)< t(A_3,B_3)$. Since $<$ is transitive, we have that $t(A_1,B_1)< t(A_3,B_3)$ and then $(A_1,B_1)<' (A_3,B_3)$.

\item Is straightforward.
\end{enumerate}
\edm

If follows almost directly by the definition of the structure in $C_s(S)$ that:

\bpr
$t:C_s(S)\to S$ is a morphism of surreal algebras.
\epr

\bre In the case of the three principal examples of SUR-algebras  we have that $t:C_s(SA)\to SA$ and $t:C_s(ST)\to ST$ are bijections and $t:C_s(No)\to No$ is a surjection.
\qdr\ere

\bpr
If $f : S \ra S'$ is a morphism then $C_s(f) : C_s(S) \ra C_s(S') : (A,B) \mapsto (f[A], f[B])$ is a morphism.
\epr

\bpr
$C_s$ determines $C_s(f)(A,B) = (f[A],f[B])$ a functor from $SUR$ to $SUR$,  and $t$ determines a natural transformation $t :Id_{SUR-alg} \to C_s$
\epr

From a direct application of the Proposition \ref{casesCut-pr}, we obtain the following:

\bpr Let $\cS = (S,<,-,*,t)$ a SUR-algebra.
\begin{enumerate}
    \item If $\cS$ is universal on the category $SUR-alg$ then the following diagram commutes:
    $$( S \overset {!}\to\to C_s(S) \overset{t}\to\to S) = (S \overset{id_S}\to\to S)$$
    \item If $\cS$ is an object of full subcategory  $SUR_T-alg \hookr SUR-alg$, of all  {\bf transitive} SUR-algebra, and is universal on the this category $SUR_T-alg$, then the following diagram commutes:
    $$( S \overset {!}\to\to C_s(S) \overset{t}\to\to S) = (S \overset{id_S}\to\to S)$$
\end{enumerate}
\epr

\bre Note that: $C_s(SA)=SA$ and $C_s(ST)=ST$.
\ere

\section{Partial Surreal Algebras and morphisms} \label{2.3-sec}

In several recursive constructions, the intermediate stages play an important role in the comprehension of the object constructed itself. As we have seen in the section 2.1, all surreal algebra is a {\em proper class} but, on the other hand, the intermediate stages of the constructions of $No, SA, ST$ are sets.
To gain some flexibility and avoid  technical difficulties, we introduce in this Section the (more general and flexible ) notion of {\em partial} surreal algebra: every SUR-algebra is a partial SUR-algebra and this new notion can be supported by a set. 
Besides simple examples, that contains in particular the intermediate stages of $No, SA, ST$, and a relativized notion of Cut (partial) SUR-algebra, we are interest on general constructions of partial SUR-algebras: for that we will consider two kinds on morphisms between them. We will perform general constructions as products, sub partial-SUR-algebra and certain  kinds of directed colimits. As an application of the latter construction, we are able to prove some universal properties satisfied by $SA$ and $ST$ (and natural generalizations), that justifies its  names of (relatively)  free SUR-algebras.

\bdf \label{pSURalg-df} A {\bf partial surreal algebra (pSUR-algebra)} is a structure ${\cal S}= (S,*,-,<,t)$ where $S$ is a class (proper or improper), 
$* \in S$, $-$ is an unary function in $S$, $<$ is a binary relation in $S$ and $t : C^{t}_s(S) \to S$ is a {\em partial function}, i.e., $C^{t}_s(S) \subseteq C_s(S)$, satisfying:

\begin{enumerate}
  \item[{\bf (pS1)}] $<$ is an acyclic relation.
  \item[{\bf (pS2)}] $\forall x\in S$, $-(-x)=x$
  \item[{\bf (pS3)}] $-*=*$.
  \item[{\bf (pS4)}] $\forall a,b\in S$, $a<b$ iff $-b<-a$
  \item[{\bf (pS5)}] If $(A,B)\in C^t_s(S)$, then $A<t(A,B)<B$.
  \item[{\bf (pS6)}] If $(A,B)\in C^t_s(S)$, then $(-B,-A)\in C^t(S)$ and $-t(A,B)=t(-B,-A)$.
	\item[{\bf (pS7)}] $(\emptyset,\emptyset) \in C^t_s(S)$ and $*=t(\emptyset,\emptyset)$.
  \end{enumerate}
	\qdr\edf

	Note that (pS1), (pS2), (pS3) and (pS4) coincide, respectively, with the SUR-algebra axioms (S1), (S2), (S3) and (S4). The statements (pS5), (pS6) and (pS7) are relative versions of, respectively, the SUR-algebra axioms (S5), (S6) and (S7). SUR-algebras are precisely the pSUR-algebras $\cS$ such that $C^{t}_s(S) = C_s(S)$.

\


\bdf \label{pSURmor-df} Let ${\cal S} = (S, <, -, *, t)$ and ${\cal S}' = (S', <', -', *', t')$ be partial SUR-algebras. Let $h : {S} \rightarrow {S}'$ be (total) function and consider the conditions below:
\begin{enumerate}
\item [{\bf (Sm1)}] $h(*) = *'$.
\item [{\bf (Sm2)}] $h(-a) = -'h(a)$, $\all a \in S$.
\item [{\bf (Sm3)}] $a<b \implies h(a) <' h(b)$, $\all a,b \in S$.
\item [{\bf (pSm4)}] $(A,B)\in C^t_s(S) \implies (h[A],h[B])\in C^{t'}_s(S')$  and  $h(t(A,B))=t'(h[A],h[B])$, $\all (A,B)\in C^t_s(S)$.
\item [{\bf (fpSm4)}]$(A,B)\in C_s(S) \implies (h[A],h[B])\in C^{t'}_s(S')$  and  $h(t(A,B))=t'(h[A],h[B])$, $\all (A,B)\in C^t_s(S)$.
\end{enumerate}

We will say that $h : {\cal S} \rightarrow {\cal S}'$ is:\\
\bu  a {\bf partial SUR-algebra morphism (pSUR-morphism)} when it satisfies: (Sm1), (Sm2), (Sm3) and (pSm4);\\
\bu  a {\bf \underline{full} partial SUR-algebra morphism (fpSUR-morphism)} when it satisfies: (Sm1), (Sm2), (Sm3) and (fpSm4).
\qdr\edf

\bre \label{pSURmor-re}

\begin{itemize}

\item Note that the property  (Sm3) entails: $(A,B)\in C_s(S) \implies (h[A],h[B])\in C_s(S')$. 

\item The conditions (Sm1), (Sm2) and (Sm3)  are already present in the definition of SUR-algebra morphism. The property:\\
\begin{enumerate}
\item [{\bf (Sm4)}] $h(t(A,B))=t'(h[A],h[B])$, $\all (A,B)\in C_s(S)$;\\
	completes the definition of SUR-algebra morphism.
	\end{enumerate}

\item Every full partial SUR-algebra morphism is  partial SUR-algebra morphism. 

\item  Let $S, S'$ be partial SUR-algebras and $h : S \to S'$  is a map. Suppose that $S$ or $S'$ is a SUR-algebra,  then $h$ is a pSUR morphism  iff $h$ is a fpSUR-morphism.

\item If $S$ is a partial SUR-algebra, then: $id_S : S \to S$ is a pSUR-morphism and $id_S : S \to S$ is a fpSUR-morphism iff $S$ is a SUR-algebra.

\item  Let $h: S \to S'$, $h' : S' \to S''$ be pSUR morphisms:\\
\bu Then $f' \circ f$ is a pSUR-morphism.\\
\bu  If $f$ is fpSUR-morphism,  then $f' \circ f$ is a fpSUR-morphism. In particular, the composition of fpSUR-morphisms is a fpSUR-morphism.

\end{itemize}
\qdr\ere

\bdf \label{CATpSUR-df} {\bf The category of partial SUR-algebras:}

We will denote by  $pSUR-alg$  the ("very-large") category such that $Obj(pSUR-alg)$ is the class of all partial SUR-algebras 
and $Mor(pSUR-alg)$ is the class of all partial SUR-algebras morphisms, endowed with obvious composition and identities.
\qdr\edf

\bre \label{CATpSUR-re}

(a) Of course, we have in the category $pSUR-alg$ the same "size issue" presented in the categories of $ZF-alg$ and $SUR-alg$: we will adopt the same "solution" explained in Remark \ref{CATSUR-re}. An alternative is to consider only "small" partial SUR-algebras (and obtain "large" category --instead of very large-- $pSUR_s-alg$, of all small partial SUR-algebras)  since we will see that there are  set-size partial SUR-algebras: we will not pursue this track because our main concern in considering partial SUR-algebras is get flexibility to make (large indexed) categorial constructions with small partial SUR-algebras to obtain a total SUR-algebra as a (co)limit process, i.e., we want $pSUR \supseteq SUR$.

(b) We saw above that, even if the class of full morphism of partial SUR-algebras is closed under composition, it does not determines a category under composition, since it lacks the identities for the small partial SUR-algebras. However this notion will be useful to perform constructions of total SUR-algebra as colimit of a large diagram small partial SUR-algebras and fpSUR-morphisms between them (see Subsection 2.3.4).
\qdr\ere

\bct \label{Sigma-str-re} Denote $\Sigma$-str the (very large) category such that:

(a) The objects of $\Sigma$-str are the structures ${\cal S}= (S,*,-,<,t)$ where $S$ is a class, $* \in S$, $-$ is an unary function in $S$, $<$ is a binary relation in $S$ and $t : D^t \to S$ is a function such that $D^t  \subseteq P_s(S)\times P_s(S)$.

(b) Let ${\cal S} = (S, <, -, *, t)$ and ${\cal S}' = (S', <', -', *', t')$ be partial SUR-algebras. A $\Sigma$-morphism, $h : {\cal S} \rightarrow {\cal S}'$,  is a (total) function $h : { S} \rightarrow { S}'$ satisfying the conditions below:
\begin{enumerate}
\item [{\bf ($\Sigma$m1)}] $h(*) = *'$.
\item [{\bf ($\Sigma$m2)}] $h(-a) = -'h(a)$, $\all a \in S$.
\item [{\bf ($\Sigma$m3)}] $a<b \implies h(a) <' h(b)$, $\all a,b \in S$.
\item [{\bf ($\Sigma$m4)}] $(h \times h)[D^t] \sub D^{t'}$ and $h(t(A,B))=t'(h[A],h[B])$, $\all (A,B)\in D^t$.
\end{enumerate}

(c) Endowed with obvious composition and identities, $\Sigma$-str is a very large category and 
$$SUR-alg \hookr pSUR-alg \hookr \Sigma-str$$
 are inclusions of full subcategories.
\qdr \ect



\subsection{Simple examples} \label{2.3.1-subsec}

In this short Subsection we just present first examples of partial SUR-algebras and its morphisms.


\bex \label{ABord-ex} Let $(G,+, -, 0, <)$ be  a linearly ordered group. For each and select $a \in G$ such that $a \geq 0$ (respect. $a \in G\cup\{\infty\}$ such that  $a>0$)then $X_a := [-a, a] \subseteq G$ (respect. $X_a := ]-a, a[ \sub G$), is a partial SUR-algebra, endowed with  obvious definitions of $*, -, <$ and such that:\\
(1) $C^t_s(X_a)  :=$  $\{(x^{<}, x^{>}) : x \in X_a\}$, $t(x^{<}, x^{>}) := x \in X_a $ ($t$ is bijective);\\
or, alternatively, \\
(2) $C^t_s(X_a)  :=$  $\{ (L,R) \in C_s(X_a) : \exists (!) x \in X_a \ L^{\leq} = x^{<}$,  $R^{\geq} = x^{>} \}$, 
$t(L,R) := x \in X_a$ ($t$ is surjective).

Note that if $b \geq a$, then the inclusion $X_a \hookr X_b$ is a $pSUR$-morphism, if $X_a, X_b$ are endowed with the second kind of $t$-map.
\qdr\eex

Another simple (and useful) class of examples are given by the ordinal steps of the recursive constructions of the SUR-algebras $SA$, $ST$ and $No$.

\

\bex  \label{pSA-ex} For any ordinal $\alpha$ we have that the $\Sigma$-structure $(SA_\alpha,<_\alpha,-_\alpha,*_\alpha,t_\alpha)$ is a partial SUR-algebra with the above definitions:
\begin{itemize}
    \item $*_\alpha=*$
    \item $-_\alpha=-\rest_{SA_\alpha}$
    \item $<_\alpha=<\rest_{SA_\alpha\times SA_\alpha}$
    \item $C_s^t(SA_\alpha)=C_s(SA^{(\alpha)})$ and $t_\alpha=t\rest_{C_s(SA^{(\alpha)})}$
\end{itemize}
\qdr \eex

Just like in the previous example, we have:

\bex \label{pST-ex}
\begin{itemize}
    \item $*_\alpha=*$
    \item $-_\alpha=-\rest_{ST_\alpha}$
    \item $<_\alpha=<\rest_{ST_\alpha\times ST_\alpha}$
    \item $C^t_s(ST_\alpha)=C_s(ST^{(\alpha)})$ and $t_\alpha=t\rest_{C_s(ST^{(\alpha)})}$
\end{itemize}
\qdr \eex

\bex \label{pNo-ex}
For any given $\alpha \in On$, the $\Sigma$-structure $(No_\alpha,*_\alpha,-_\alpha,<_\alpha,t_\alpha)$ is a partial SUR-algebra with the operations defined above:
\begin{itemize}
    \item $*_\alpha=0$
    \item $-_\alpha=-\rest_{No_\alpha}$
    \item $<_\alpha=<\rest_{No_\alpha\times No_\alpha}$
    \item $C^t_s(No_\alpha)=C_s(No^{(\alpha)})$ and $t_\alpha=t\rest_{C_s(No^{(\alpha)})}$
\end{itemize}
\qdr \eex

\bre \label{fpSUR-re}

\bu Note that in the three examples above $S = SA, ST, No$, the inclusion $S_\alpha \hookr S_\beta$ is a pSUR-morphism, where $\alpha \leq \beta \leq \infty$ are "extended" ordinals, with the convention $S_\infty := S$.

\bu We can also define partial SUR-algebras on the sets $SA^{(\alpha)}, ST^{(\alpha)}, No^{(\alpha)}$, for each $\alpha \in On\setminus\{0\}$ 
(this is useful!).

\bu Note that $i_\alpha: SA^{(\alpha)} \hookr SA_\alpha$ is a fpSUR-algebra morphism, for each $\alpha \in On\setminus\{0\}$. It can be established, by induction on $\alpha \in On\setminus\{0\}$ that for each $\gamma < \alpha$ $i_{\gamma\alpha} : SA_\gamma \hookr SA_\alpha$ is a fpSUR-morphism.  An analogous situation occurs to the partial SUR-algebras $ST_\gamma \hookr ST^{(\alpha)} \hookr ST_\alpha$.
\qdr\ere



\subsection{Cut partial Surreal Algebras} \label{2.3.2-subsec}

In this short Subsection we present an adaption/generalization  of the notion of "Cut Surreal Algebra", introduced in the Subsection 2.2.4, to the realm of {\em partial} SUR-algebra.

\bdf \label{CUTpSUR-df} Let $\cS = (S,<,-,*,t)$ be a partial SUR-algebra. The {\em Cut structure} of $\cS$ is the $\Sigma$-structure  ${\cS}^{(t)} = (S', <',-',<', t')$, where:

\begin{enumerate}
\item $S' := C_s^t(S)$
    \item $*':= (\emptyset,\emptyset)$
    \item $-'(A,B) := (-B,-A)$
    \item $(A,B)<'(C,D) \iff t(A,B) < t(C,D)$
		\item $\all \alpha, \beta \sub C^t_s(S)$, $(\alpha, \beta) \in dom(t')$ iff $\alpha<'\beta$ and $(t[\alpha], t[\beta]) \in dom(t)$
		\item $t' : C^{t'}_s(C^t_s(S)) \to C^t_s(S)$, $(\alpha,\beta) \mapsto t'(\alpha,\beta) := (t[\alpha], t[\beta])$
\end{enumerate}
\qdr\edf

The list below a sequence of results on Cut Partial SUR-algebras that extend the results presented in the Subsection 2.2.4 on Cut SUR-algebras: its  proofs will be omitted.

\bpr \label{CutpSUR-pr} Let $\cS = (S,<,-,*,t)$ be a partial SUR-algebra. Then:

(a) ${\cS}^{(t)} = (S', <',-',<', t')$ as defined above is a partial SUR-algebra. Moreover, if $\cS$ is a SUR-algebra, i.e. $C_s^t(S) = C_s(S)$, then  ${\cS}^{(t)}$ is a SUR-algebra, i.e. $C^{t'}_s(C_s^t(S)) = C_s(C_s(S))$.

(b)  $t : C_s^t(S) \to S$ is a morphism of partial SUR-algebras. Moreover, if $\cS$ is a SUR-algebra, then $t$ is a fpSUR-algebra morphism.
\qdr\epr

\bpr \label{pSURkind-pr} Let $\cS = (S,<,-,*,t)$ be a partial SUR-algebra. Then:

(a) If $S$ is transitive, then $C^t_s(S)$ is transitive.
    
(b) If $S$ is linear, then $C^t_s(S)$ is pre-linear\footnote{I.e., denote $\sim_t$ the equivalence relation on $C^t_s(S)$ given by $(A,B) \sim_t (C,D)$ iff $t(A,B) = t(C,D)$. Then it holds exactly one between of the alternatives: $(A,B) <' (C,D)$; $(A,B) \sim_t (C,D)$; $(C,D) <' (A,B)$.}.
\qdr\epr

\bpr \label{pSURfunct-pr}

(a) If $f : \cS \ra {\cS}'$ is a morphism of partial SUR-algebras then $C_s^t(f) : C^t_s(S) \ra C^t_s(S')$, given by: $(A,B) \mapsto (f[A], f[B])$ is a morphism of partial SUR-algebras.

(b) The cut partial SUR-algebra construction determines a (covariant) functor $C^t_s : pSUR \to pSUR$:
$$ (S \overset{f}\to\to S') \ \mapsto \  (C^t_s(S) \overset{C_s^t(f)}\to\to  C^t_s(S'))$$

(c) The $t$-map  determines a natural transformation between functors on $pSUR-alg$,  $t :Id_{pSUR-alg} \to C^t_s$.
\qdr\epr



\subsection{Simple constructions on pSUR}  \label{2.3.3-subsec}

In this Section, we will verify the full subcategory $pSUR-alg \hookr \Sigma-str$ is closed under some simple categorial constructions: as  ($\Sigma$-)substructure and non-empty products. We also present some results on initial objects and (weakly) terminal objects.

We can also define a notion of substructure in the category $pSUR$:

\bdf \label{Sigmasub-df} Let ${\cal S}=(S,<,-,*,t)$ and ${\cal S}'=(S',<',-',*',t')$ be $\Sigma$-structures. $\cS$ will be called a $\Sigma$-substructure of $\cS$ whenever:\\
(s1) $S\subseteq S'$;\\
(s2) $<=<'_{\rest_{S\times S}}$;\\
(s3) $-=-'_{\rest_{S\times S}}$;\\
(s4) $*=*'$;\\
(s5) $dom(t) = t'^{-1}[S] \cap (P_s(S) \times P_s(S)) := \{(A,B) \in dom(t') \cap (P_s(S) \times P_s(S)): t'(A,B) \in S\} \sub dom(t')$ and $t=t'_{\rest} : \dom(t) \to S$.
  \qdr\edf

\bre \label{pSURsub-re}

(a) The inclusion $i : S \hookr S'$ determines a $\Sigma$-morphism.

(b) By conditions (s1) and (s2) above note that $C_s(S,<) = C_s(S',<') \cap (P_s(S) \times P_s(S))$.

(c)  By item (b): if $dom(t') \sub C_s(S',<')$, then $dom(t) \sub C_s(S,<)$. 

(d) By the results presented in the Subsections 2.2.2 and 2.2.3, for any two extends ordinals $\alpha\leq \beta \leq \infty$ we have:\\
\bu $SA_\alpha$ is a $\Sigma$-substructure of $SA_\beta$.\\
\bu $ST_\alpha$ is a $\Sigma$-substructure of $ST_\beta$.

(e) An useful generalization of the notion of $\Sigma$-substructure is the notion of $\Sigma$-embedding: a $\Sigma$-morphism $j : \cS \to \cS'$ is a $\Sigma$-embedding when: \\
(e1) it is injective;\\
(e2) $\all a,b \in S, (a < b \Lra j(a)<'j(b))$;\\
(e3) $\all (A,B) \in P_s(S) \times P_s(S), ( (A,B) \in dom(t) \Lra t'(j[A],j[B]) \in range(j) )$. 

(e) An inclusion $i : S \hookr S'$ determines a $\Sigma$-embedding precisely when $\cS$ is a $\Sigma$-substructure of $\cS'$.  Note that the $\Sigma$-embeddings $j : \cS \to \cS'$ are precisely the $\Sigma$-morphisms  described (uniquely) as $j = i \circ h$, where $i : \cS^j \hookr \cS'$ is a $\Sigma$-substructure inclusion and $h : \cS \to \cS^j$ is a $\Sigma$-isomorphism.

(f) For technical reasons, we consider an even more general notion: a $\Sigma$-morphism  $j : \cS \to \cS'$ is a
$\Sigma-quasi$-embedding  whenever it satisfies the conditions (e1) and (e3) above.
\qdr\ere

By a straightforward verification we obtain the:

\bpr \label{pSURsub-pr} Let  $j : \cS \to \cS'$ be a $\Sigma$-embedding of $\Sigma$-structures. If ${\cal S}'$ is a partial SUR-algebra, then ${\cal S}$ is a partial SUR-algebra.  
\qdr\epr




\bdf \label{pSURprod-df} Given a {\em non-empty} indexed {\em set} of partial $\Sigma$-structure $\cS_i=(S_i,<_i,-_i,*_i,t_i)$, ${i\in I}$, we define the $\Sigma$-structure product  $\cS = (S,<,-,*, t)$ as follows: Let ${\cal S}=(S,<,-,*,t)$ and ${\cal S}'=(S',<',-',*',t')$ be $\Sigma$-structures. $\cS$ will be called a $\Sigma$-substructure of $\cS$ whenever:\\
(a) $S = \prod_{i \in I} S_i$;\\
(b) $< = \{((a_i)_{i \in I}, (b_i)_{i \in I}):  a_i <_i b_i, \all i \in I\}$;\\
(c) $-(a_i)_{i \in I} = (-_i a_i)_{i \in I}$;\\
(d) $* = (*_i)_{i \in I}$;\\
(e) $dom(t) = \bigcap_{i \in I}(\pi_i \times \pi_i)^{-1}[dom(t_i)]  = \{((A_i){i \in I},(B_i)_{i \in I}) \in  P_s(S) \times P_s(S)): (A_i, B_i) \in dom(t_i), \all i \in I\}$ and $t((A_i){i \in I},(B_i)_{i \in I})) = (t_i(A_i,B_i))_{i \in I}$.
\qdr\edf

Note that: For each $i \in I$, the projection $\pi_i : S \to S_i$ is a $\Sigma$-structure morphism.

By a straightforward verification we obtain:

\bpr \label{pSURprod-pr} Keeping the notation above.

(a) The pair $(\cS, (\pi)_{i \in I})$ above defined constitutes a(the) categorial product in $\Sigma$-str. I.e., for each diagram $({\cS}', (f_i)_{i \in I})$ in $\Sigma$-str such that $f_i : {\cS}' \to \cS_i$, $\all i \in I$, there is a unique $\Sigma$-morphism $f : \cS' \to \cS$ such that $\pi_i \circ f = f_i, \all i \in I$.  

(b) Suppose that $\{\cS_i : i \in I\} \sub$ pSUR-alg. Then $\cS \in$ pSUR-alg and  $(\cS, (\pi)_{i \in I})$ is the product in the category pSUR-alg. 
\qdr\epr

\bpr \label{linpSUR-pr} Let $f : \cS \to {\cS}'$ be a pSUR-alg morphism. If $(S,<)$ is  strictly linearly ordered, then:\\
(a)  $\all a, b \in S$, $a<b \iff f(a)<'f(b)$; \\
(b) $f$ is an injective function.
\epr
\bdm
If $a<b$, then $f(a)<'f(b)$, since $f$ is a $Sigma$-structure morphism. Suppose that $f(a)<'f(b)$ but $a \not < b$, then $a = b$ or $b<a$, thus $f(a) = f(b)$ or $f(b) <' f(a)$. In the case, we get a contradiction with $f(a)<'f(b)$, since $<'$ is an acyclic relation. This establishes item (a). Item (b) is similar, since $<$ satisfies trichotomy and $<'$ is acyclic.
\qdr\edm


The result above yields some information  concerning the empty product (= terminal object) in pSUR-algebras.

\bpr \label{termpSUR-co}
If there exists a {\em weakly} terminal object\footnote{Recall that an object in a category is weakly terminal when it is the target of {\em some} arrow departing from each object of the category.} $\cS_1$ in the category pSUR-alg then $\cS_1$ must be a proper class.
\epr
\bdm Suppose that $\cS_1$ is an weakly terminal object in pSUR-alg. Since the (proper class) SUR-algebra $No$ is strictly linearly, then by Proposition \ref{linpSUR-pr} above anyone of the existing morphisms $f : No \to \cS_1$ is injective. Then $\cS_1$ (and $C^t_s(\cS_1)$) must be a proper class. 
\qdr\edm

If we consider the small size version of pSUR, we can guarantee by an another application of Proposition \ref{linpSUR-pr}, that this (large but not very-large) category does not have  (weakly) terminal objects: there are small abelian linearly ordered abelian groups (or even the additive part of a ordered/real closed field) of arbitrary large cardinality, and we have seen in Example \ref{ABord-ex} how to produce small pSUR-algebras from that structures.

\

Concerning initial objects we have the following:

\bpr \label{inipSUR-pr} 

(a) Consider the $\Sigma$-structure $\cS_0 = (S_0, *, -, <, t)$ over a singleton set $S_0 := \{*\}$, with $< := \empty$,  $D^t = dom(t) := \empty$ (thus $\cS_0 \notin pSUR-alg$) and with $- : S_0 \to S_0$ and $t : D^t \to S_0$ the unique functions available. Then $\cS_0$ is the (unique up to unique isomorphism) initial object in $\Sigma$-str.

(b) Consider the $\Sigma$-structure $\cS^p_0 = (S_0, *, -, <, t^p)$  over a singleton set $S_0 := \{*\}$, with $< := \empty$,  $D^{t^p} = dom(t^p) := \{(\empty,\empty)\} \sub C_s(S_0,<)$ and with $- : S_0 \to S_0$ and $t^p : D^{t^p} \to S_0$ the unique functions available. Then $\cS^p_0$ is the (unique up to unique isomorphism) initial object in pSUR-alg. 

\epr
\bdm 

(a) Let ${\cal S}'$ be a $\Sigma$-structure and let $h: \{*\} \to S'$ be the unique function such that $h(*) = *' \in S'$, then clearly $h$ is the unique $\Sigma$-structure morphism from $\cS_0$ into ${\cal S}'$: note that $(h \times h)_\rest : dom(t) = \empty \to dom (t')$ is such that $t' \circ (h \times h)_\rest = h \circ t$.

(b) It is easy to see that $\cS^p_0$ is a partial SUR-algebra. Let ${\cal S}'$ be a partial SUR-algebra and let $h: \{*\} \to S'$ be the unique function such that $h(*) = *' \in S'$, then clearly $h$ is the unique $\Sigma$-structure morphism from $\cS_0$ into ${\cal S}'$: since $(\empty,\empty) \in dom(t')$, note that $(h \times h)_\rest : dom(t^p) := \{(\empty,\empty)\}  \to dom(t')$  is such that $t' \circ (h \times h)_\rest = h \circ t$.
\qdr\edm



%



\subsection{Directed colimits of partial Surreal Algebras} \label{Dircolim-subsec} 

One of the main general constructions in Mathematics is the colimit of an upward directed diagram. In the realm of partial SUR-algebras this turns out to be essential for the constructions of SUR-algebras and to obtain general results about them. We can recognize the utility of this process by the cumulative constructions of our main examples: $No, SA, ST$. Thus we will be concerned only with the  colimit of {\bf small} partial SUR-algebras, but over a possibly a {\em large} directed diagram. 

This Section is completely technical: we intend provide full  proofs and register them  possibly in a Appendix. Its consequences/applications  are very interesting: see the entire Section 2.4 and the Theorems at the end of Section 3.1.

Recall that:

\bu Given a regular "extended" cardinal $\kappa$ (where a "$card(X) = \infty$" means that $X$ is a proper class\footnote{Recall that in NBG, all the proper classes are in bijection, by the global form of the axiom of choice.}), a partially ordered class  $(I, \leq)$ will be $\kappa$-directed, if every subclass $I' \sub I$ such that $card(I') < \kappa$ admits an upper bound in $I$.  be a $\omega$-directed ordered class.

\bu $pSUR_s-alg$ denotes the full subcategory of $pSUR-alg$  determined by  of all {\em small} partial SUR-algebras and its morphisms (then $SUR-alg \cap pSUR_s-alg = \empty$). Analogously, we will denote $\Sigma_s$-str the full subcategory of $\Sigma$-str  determined by  of all {\em small} partial $\Sigma$-structures and its morphisms.

\

\bct \label{colim-ct} {\bf The (first-order) directed colimit construction:} Let  $(I, \leq)$ is a $\omega$-directed ordered class and consider ${\cal D} : (I,\leq) \to \Sigma_s-str$, $(i \leq j) \mapsto (\cS_i \overset{h_{ij}}\to\to \cS_j)$ be a diagram. Define: \\
\bu $S_\infty := (\sqcup_{i \in I} S_i)/ \equiv$, the set-theoretical colimit, i.e. $\equiv$ is the least equivalence relation on the class $\sqcup_{i \in I} S_i$ such that $(a_i,i) \equiv (a_j,j)$ iff there is $k \geq i,j$ such that $h_{ik}(a_i) = h_{jk}(a_j) \in S_k$;\\
\bu $h_j : S_j \to S_\infty$, $a_j \mapsto [(a_j,j)]$;\\
\bu $* := [(*_i,i)]$ (= $[(*_j, j)]$, $\all i,j \in I$);\\
\bu $-[(a_i,i)] = [(-_i a_i,i)]$;\\
\bu $[(a_i,i)] < [(a_j,j)]$ iff there is $k \geq i,j$ such that $h_{ik}(a_i) <_k h_{jk}(a_j) \in S_k$
 \qdr\ect


With the construction above, it is straightforward to verify that $(S_\infty, <, -, *)$ is the colimit in the appropriate category of {\em first-order} (but possibly large) structures\footnote{I.e., we drop the second-order part of the $\Sigma$-structure: the map $t : dom(t) \sub P_s(S_\infty)\times P_s(S_\infty) \to S\infty$.}, with colimit co-cone  $(h_j : S_j \to S_\infty)_{j \in I}$ and, if ${\cal D} : (I,\leq) \to pSUR_s-alg$, then the same (colimit) co-cone is in the "first-order part" of the category $pSUR-alg$, i.e., it satisfies the properties [pS1]--[pS4] presented in Definition \ref{pSURalg-df}. However, to "complete" the $\Sigma$-structure (respect. pSUR-algebra)  we will need some extra conditions below: 
d


\bpr \label{Sigmacolim-pr} Let ${\cal D} : (I,\leq) \to \Sigma_s-str$, $(i \leq j) \mapsto (\cS_i \overset{h_{ij}}\to\to \cS_j)$ be a diagram such that:\\
(i)  $(I, \leq)$ is a $\omega$-directed ordered class and $h_{ij} : S_i \to S_j$ is a injective $\Sigma$-morphism, whenever $i\leq j$;\\
 or;\\
(ii)  $(I, \leq)$ is a $\infty$-directed ordered class (e.g. $(On, \leq)$), \\
then  $S_\infty := (\sqcup_{i \in I} S_i)/ \equiv$ is a (possibly large) partial $\Sigma$-structure and $(h_j : S_j \to S_\infty)_{j \in I}$ is a colimit cone in the category $\Sigma-str$.
\qdr\epr

\bpr \label{pSURcolim-pr} If ${\cal D} : (I,\leq) \to pSUR_s-alg$, $(i \leq j) \mapsto (\cS_i \overset{h_{ij}}\to\to \cS_j)$  is a diagram, where:\\
(i)  $(I, \leq)$ is a $\omega$-directed ordered class and $h_{ij} : S_i \to S_j$ is a injective pSUR-morphism, whenever $i\leq j$;\\
 or;\\
(ii)  $(I, \leq)$ is a $\infty$-directed ordered class (e.g. $(On, \leq)$) \\
then  $S_\infty$ is a (possibly large) partial SUR-algebra and $(h_j : S_j \to S_\infty)_{j \in I}$ is a colimit cone in the category $pSUR-alg$.
\qdr\epr

\bpr \label{fpSURcolim-pr} The subclass of morphisms $fpSUR \subseteq pSUR$ is {\em closed} under directed colimits in the cases (i) and (ii) described in the Proposition above. More precisely: if ${\cal D} : (I,\leq) \to pSUR_s-alg$ is a directed diagram  satisfying (i) or/and (ii)above and such that $h_{ij} : S_i \to S_j$ is a fpSUR-morphism, whenever {\underline{$i < j$}},  then the colimit 
co-cone $\all j \in I$, $(h_j : S_j \to S_\infty)_{j \in I}$ is a formed by  fpSUR-algebra morphisms. Moreover:\\
(a) If $(I,\leq)$ is $\infty$-directed and the transition arrows $(h_{ij})_{i \leq j}$ are injective, then $S_\infty$ is a SUR-algebra (thus it is a proper class);\\
(b) If the transition arrows $(h_{ij})_{i \leq j}$ are injective (respect. $\Sigma-quasi$-embedding, $\Sigma$-embedding), then the cocone arrows  $(h_{j})_{j \in I}$ are injective (respect. $\Sigma-quasi$-embedding, $\Sigma$-embedding);\\
(c) If $t_i : C^{t_i}_s(S_i) \to S_i$ is injective (respect. surjective/bijective), $\all i \in I$, then $t^\infty : C^{t^\infty}_s(S^\infty) \to S^\infty$ is injective (respect. surjective/bijective).
\qdr\epr

\bex \label{colim-ex}

We have noted in Remark \ref{fpSUR-re} that for each sequence of ordinal $\gamma < \beta < \alpha$, $i_{\gamma\beta} : SA_\gamma \hookr SA_\beta$ is fpSUR-algebra morphism. It is also a $\Sigma$-embedding.  Then, for each $\alpha >0$, $SA^{(\alpha)} \cong colim_{\gamma < \alpha} SA_\gamma$ as a  pSUR-algebra and $i^{(\alpha)}_\gamma : SA_\gamma \hookr SA^{(\alpha)}$ determines a colimit co-cone of an $\omega$-directed diagram\footnote{In fact it is $\kappa$ directed diagram, where $\kappa$ is any regular cardinal such that $\kappa \leq \alpha+\omega$.} formed by  fpSUR-algebras embeddings. Moreover $SA = SA_\infty \cong colim_{\gamma \in On} SA_\gamma$ and $i^{\infty}_\gamma : SA_\gamma \hookr SA$ determines a colimit co-cone a $\infty$-directed diagram  over formed by  fpSUR-algebras embeddings.

Analogous results holds for $ST^{(\alpha)} \cong colim_{\gamma < \alpha} ST_\gamma$, $\alpha >0$,  and $ST \cong colim_{\gamma \in On} ST_\gamma$.
\qdr\eex



\section{Universal Surreal Algebras}  \label{sec2.4}


In this  section, we present some {\em categorical-theoretic} universal properties\footnote{An analysis of model-theoretic universal properties of the "first-order part" of (partial) SUR-algebras, and its possible connections with categorial-theoretic universality presented here,  will be theme of future research, see section \ref{conclusion} for more details.} concerning  SUR-algebras and partial SUR-algebras. We will need notions, constructions and results developed in the previous sections to provide, for each small {\em partial} SUR-algebra $I$, a "best" SUR-algebra over $I$, $SA(I)$, (respect. a "best" {\em transitive} SUR-algebra over $I$, $ST(I)$).  As a consequence of this result (and its proof) we will determine the SUR-algebras $SA$ and $ST$ in the category of SUR by universal properties that characterizes them uniquely up to unique isomorphisms: these will justify the adopted names "$SA$ = the free surreal algebra" and "$ST$ = the free transitive surreal algebra".

We start with the following

\bct \label{IpSUR-ct} {\bf Main construction:} Let $\cI = (I, *,-,<, t)$ be a  partial SUR-algebra. Consider:

(a) The set-theoretical pushout diagram over $(I \overset{t}\to{\leftarrow} C_s^t(I) \overset{incl}\to{\hookr} C_s(I))$:

\begin{picture}(140,140)

\setlength{\unitlength}{.7\unitlength}
\thicklines
\put(12,20){$C^t_s(I)$}
\put(12,162){$C_s(I)$}
\put(177,20){$I$}
\put(157,162){$(I \sqcup C_s(I))/\sim$}
 \put(65,20){\vector(1,0){100}}
\put(65,170){\vector(1,0){80}}
\put(185,40){\vector(0,1){110}}
\put(45,40){\vector(0,1){110}}
\put(5,90){$incl$}
\put(200,90){$i_0$}
\put(110,0){$t$}
\put(110,180){$i_1$}
\end{picture}

Note that:

\bu $(I \sqcup C_s(I))/\sim$, is the vertex of the set-theoretical pushout diagram, where $\sim$ is the least equivalence relation\footnote{Recall that the least equivalence relation on a set $X$ that contains $R \sub X \times X$ is obtained from $R$ adding the opposite relation $R^{-1}$ and the diagonal relation $\Delta_X$, and then taking the transitive closure  $trcl(R \cup R^{-1} \cup \Delta_X) = R^{(eq,X)}$.}   
on $I \sqcup C_s(I)$ such that $(x,0) \sim ((A,B),1)$ iff $(A,B) \in C^t_s(I)$ and $x = t(A,B)$. \\
\bu If $I$ is small, then $I^+$ is small.\\
\bu $\all (A,B), (C,D)  \in  C_s(I) \setminus C^t_s(I)$, $((A,B),1) \sim ((C,D),1)$ iff $(A,B) = (C,D)$ (by induction on the number of steps that witness the transitive closure).\\
\bu $\all x, y  \in I$, $(x,0) \sim (y,0)$ iff $x = y$ (by induction on the number of steps needed in the transitive closure).\\
\bu Since $C^t_s(I) \hookr C_s(I)$ is injective function, then $i_0 : I \to (I \sqcup C_s(I))/\sim$, $x \mapsto [(x,0)]$ is an injective function (see above) and $(i_0)^+ : (P_s(I) \times P_s(I)) \mapsto (A,B) \mapsto (i_0[A],i_0[B])$ is an injective function.\\

(b)  ${\cI}^+ := (I^+, *^+,-^+,<^+, t^+)$ the $\Sigma$-structure defined below:\\
\bu $I^+ := (I \sqcup C_s(I))/\sim$.\\
\bu $*^+ := [(*,0)] = [((\empty, \empty),1)]$.\\
\bu $-^+[(x,0)] := [(-x,0)]$; \\
$-^+[((A,B),1)] := [((-B,-A),1)]$.\\
\bu Define $<^+$ by cases (only three):\\
$[(x,0)] <^+ [(y,0)]$ iff $x<y$;\\
$[(x,0)] <^+ [((A,B),1)]$ iff $x \in A$, whenever $(A,B) \in  C_s(I) \setminus C^t_s(I)$;\\
$[((A,B),1)] <^+ [(y,0)]$ iff $y \in B$,  whenever $(A,B) \in  C_s(I) \setminus C^t_s(I)$.\\  
 Note that $C_s(i_0) = {(i_0)^+}_\rest : C_s(I) \to C_s(I^+)$, $(A,B) \mapsto (i_0[A], i_0[B])$, is an injective function with adequate domain and codomain.\\
\bu Define $C^{t^+}_s(I^+) := range(C_s(i_0)) \sub C_s(I^+)$ (thus $C_s(I) \cong C_s^{t^+}(I^+)$) and $t^+ : C^{t^+}_s(I^+) \to I^+$, $(i_0[A], i_0[B]) \mapsto t^+(i_0[A], i_0[B]) := [((A,B),1)]$.\\
Note that $(A,B) \in dom(t)$ iff $t^+(i_0[A], i_0[B]) \in range(i_0)$.

\

Thus we obtain another set-theoretical pushout diagram that is isomorphic to the previous pushout diagram:

\begin{picture}(140,140)

\setlength{\unitlength}{.7\unitlength}
\thicklines
\put(12,20){$C^t_s(I)$}
\put(-02,162){$C_s^{t^+}(I^{+})$}
\put(177,20){$I$}
\put(177,162){$I^+$}
 \put(65,20){\vector(1,0){100}}
\put(185,40){\vector(0,1){110}}
\put(65,170){\vector(1,0){100}}
\put(45,40){\vector(0,1){110}}
\put(-25,90){$(i_0 \times i_0)_\rest$}
\put(200,90){$i_0$}
\put(110,0){$t$}
\put(110,180){$t^+$}
\end{picture}
\qdr\ect

We describe below the main technical result in this Section:

\ble \label{IpSUR-le}  Let $\cI = (I, *,-,<, t)$ be a (small) partial SUR-algebra and keep the notation in \ref{IpSUR-ct}above. Then

(a) ${\cI}^+ = (I^+, *^+,-^+,<^+, t^+)$ is a (small) partial SUR-algebra. 

(b) $i_0 : I \to I^+$ is a $\Sigma$-embedding and full morphism of partial SUR-algebras.

(c) If $t : C^t_s(I) \to I$ is injective (respect. surjective/bijective), then $t^+ : C^{t^+}_s(I^+) \to I^+$ is injective (respect. surjective/bijective).

(d) If ${\cS}' = (S',*',-',<',t') $ is a partial SUR-algebra, then for each fpSUR-algebra morphism $f : \cI \to {\cS}'$ there is a unique pSUR-algebra morphism $f^+ : {\cI}^+ \to {\cS}'$ such that $f^+ \circ i_0 = f$. In particular, if $\cS$ is a SUR-algebra, then $f$ and $f^+$ are automatically fpSUR-algebras morphisms. Moreover:\\
\bu If $t'$ is injective, then  $f$ is a $\Sigma-quasi$-embedding iff $f^+$ is a $\Sigma-quasi$-embedding.

\ele

\bdm Items (a), (b) and (c) are straightforward verifications. We will just sketch the proof of the universal property in item (d). \\

{\underline{Candidate and uniqueness:}} \\
Suppose that there is a pSUR-algebra morphism $f^+ : {\cI}^+ \to {\cS}'$ such that $f^+ \circ i_0 = f$. Since $f : I \to S$ is a full partial SUR-algebra morphism, we have $(f \times f)_\rest : C_s(I) \to C_s^{t'}(S')$. Then $(f^+ \times f^+)_\rest : C_s^{t^+}(I^+) \to C_s^{t'}(S')$ :  $(\Gamma, \Delta)= (i_0[A], i_0[B]) \mapsto (f^+[\Gamma], f^+[\Delta]) = (f[A],f[B]) \in C_s^{t'}(S')$ and $f^+(t^+(\Gamma, \Delta)) = t'(f^+[\Gamma], f^+[\Delta]) =$\\
$t'((f[A],f[B])) \in S'$. Since $range(i_0) \cup range(i_1) = S^+$,  the function $f^+$ is determined by $f$:\\
\bu $f^+(z) = f(x) \in S'$, whenever $z = [(x,0)] \in range(i_0)$;\\
\bu $f^+(z) = t'((f[A],f[B])) \in S'$, whenever $z = ([(A,B),1)] \in range(i_1)$.

\underline{Existence:}\\
Since $f : I \to S$ is a full partial SUR-algebra morphism, we have $(f \times f)_\rest : C_s(I) \to C_s^{t'}(S')$, then the arrows
$$(C_s(I) \overset{t' \circ (f \times f)_\rest}\to{\longrightarrow} S' \overset{f}\to{\longleftarrow} I)$$
yields a commutative co-cone over the diagram
$$(I \overset{t}\to{\leftarrow} C_s^t(I) \overset{incl}\to{\hookr} C_s(I)).$$
By the universal property of set-theoretical pushout, there is a unique function $f^+ : I^+ \to S'$ such that:\\
\bu  $f^+ \circ i_0 = f$;\\
\bu $f^+ \circ i_1 = t' \circ (f \times f)_\rest$.\\
Thus it remains only to check that $f' : I^+ \to S'$ is a pSUR-algebra morphism:\\
\bu $f^+(*^+) = f^+([*,0]) = f(*) = *'$;\\
\bu $f^+(-^+[(x,0)]) = f(-x) = -' f(x) = -'f^+([(x,0)])$;\\
$f^+(-^+[((A,B),1)]) = t'((f \times f)_\rest(-B,-A)) = t'(f[-B],f[-A]) = -'t'(f[B],f[A]) = -'f^+([((A,B),1)])$.\\
\bu If $[(x,0)] <^+ [(y,0)]$, then $x<y$ thus $f^+([(x,0)]) = f(x) <' f(y) = f^+([(y,0)])$;\\ 
If $(A,B) \in  C_s(I) \setminus C^t_s(I)$:\\
- if $[(x,0)] <^+ [((A,B),1)]$, then $x \in A$ and $f(x) \in f[A]$. Since$(f[A], f[B]) \in C^{t'}_s(S')$, thus $f^+([(x,0)]) = f(x) <' t'(f[A],f[B]) = f^+([((A,B),1)])$;\\
- if $[((A,B),1)] <^+ [(y,0)]$,  then $y \in B$ and $f(y) \in f[B]$. Since  $(f[A], f[B]) \in C^{t'}_s(S')$, thus $f^+([((A,B),1)]) = t'(f[A],f[B]) <' f(y) = f^+([(y,0)])$.\\
\bu If $(\Gamma, \Delta)= (i_0[A], i_0[B]) \in C^{t^+}_s(I^+)$, then $(f^+[\Gamma], f^+[\Delta]) = (f[A],f[B]) \in C^{t'}_s(S')$ and 
$f^+(t^+(\Gamma, \Delta)) = f^+(t^+((i_0[A], i_0[B])) = f^+([((A,B),1)]) = t' \circ (f \times f)_\rest([((A,B),1)]) = t'(f[A],f[B]) = t'(f^+[i_0[A]], f^+[i_0[B]]) = t'(f^+[\Gamma], f^+[\Delta])$.
\qdr\edm

\bre \label{NOM-re} 

In  the setting above, we can interpret the Conway's notions  in a very natural way:\\
\bu $Old(I) := i_0[I]  \cong I$; \\
\bu $Made(I) := I^+$;\\
\bu  $New(I):= I^+ \setminus i_0[I] = New(I)$.

Note that if $t : C^t_s(I) \to I$ is surjective (e.g. $I = No^{(\alpha)},SA^{(\alpha)},ST^{(\alpha)}$, $\alpha \in On\setminus\{0\}$), then $t^+ : C^{t^+}_s(I^+) \to I^+$ and every "made member" is  represented buy a Conway cut in of "old members". This representation is unique, whenever $t : C^t_s(I) \to I$ is bijective (e.g., $I = SA^{(\alpha)},ST^{(\alpha)}$, $\alpha \in On\setminus\{0\}$).

When $I = SA^{(\alpha)}, \alpha \in On\setminus\{0\}$ and 
$C_s^t(I) = \{(A,B) \in C_s(SA^{(\alpha)}, <_{(\alpha)}) : t(A,B) = \ll A, B \rr \in SA^{(\alpha)}$
($t : C_s(I) \to I$ is bijective), then $t^+ : C^{t^+}_s(I^+) \to I^+$ can be identified with the (bijective) map $C_s(SA^{(\alpha)}, <_{(\alpha)}) \to SA_\alpha$. 
\qdr\ere

A slight  modification in the construction of the $\Sigma$-structure presented in \ref{IpSUR-ct} above, just replacing $<^+$ by $<^{+}_{(tc)}:= trcl(<^+)$, yields the following:

\ble \label{IpSURtr-le}   Let $\cI = (I, *,-,<, t)$ be a (small) partial SUR-algebra and keep the notation in \ref{IpSUR-ct} above. Then

(a) ${\cI}^+_{(tc)} = (I^+, *^+,-^+,<^{+}_{(tc)}, t^+)$ is a (small) {\underline{transitive}} partial SUR-algebra. 

(b) $i_0 : I \to I^+_{(tc)}$ is a  $\Sigma-quasi$-embedding (see Remark \ref{pSURsub-re}.(f))  and  full morphism of partial SUR-algebras. Moreover, if $\cI$ is a {\underline{transitive}} SUR-algebra, then $i_0 : I \to I^+_{(tc)}$ is a $\Sigma$-embedding.

(c) If $t : C^t_s(I) \to I$ is injective (respect. surjective/bijective), then $t^+ : C^{t^+}_s(I^+) \to I^+$ is injective (respect. surjective/bijective).

(d) If ${\cS}' = (S',*',-',<',t') $ is a partial {\underline{transitive}} SUR-algebra, then for each fpSUR-algebra morphism $f : \cI \to {\cS}'$ there is a unique pSUR-algebra morphism $f^+ : {\cI}^+ \to {\cS}'$ such that $f^+ \circ i_0 = f$. In particular, if $\cS$ is a {{transitive}} SUR-algebra, then $f$ and $f^+$ are automatically fpSUR-algebras morphisms. 
\qdr\ele

When $I = ST^{(\alpha)}, \alpha \in On\setminus\{0\}$ and 
$C_s^t(I) = \{(A,B) \in C_s(ST^{(\alpha)}, <_{(\alpha)}) : t(A,B) = \ll A, B \rr \in ST^{(\alpha)}$
($t : C_s(I) \to I$ is bijective), then $t^+ : C^{t^+}_s(I^+_{(tc)}) \to I^+_{(tc)}$ can be identified with the (bijective) map $C_s(ST^{(\alpha)}, <_{(\alpha)}) \to ST_\alpha$.

\

\bre \label{SASTstable-re}

Note that applying the construction $( \ )^{+}$ to the SUR-algebra $SA$ we obtain $(SA)^+ \cong C_s(SA) = SA$.

Applying both constructions $( \ )^{+}$ and  $( \ )^+_{(tc)}$ to the SUR-algebra $ST$ we obtain $(ST)^+ = (ST)^+_{(tc)} \cong C_s(ST) = ST$. 
\qdr \ere

Now we are ready to state and prove the main result of this Section:

\bte  \label{adj-te} Let $\cI$ be any  small partial SUR-algebra. Then there exists SUR-algebras denoted by $SA(\cI)$ and $ST(\cI)$, and pSUR-morphisms $j^A_I : \cI \to SA(\cI)$ and $j^T_I : \cI \to ST(\cI)$ such that:

(a)  \\
(a1) $j^A_I$ is a fpSUR-morphism and a $\Sigma$-embedding;\\
(a2) If $t : C^{t}_s(I) \to I$ is injective (respect. surjective/bijective), then $t^\infty : C^{t^\infty}_s(SA(I)) \to SA(I)$ is injective (respect. surjective/bijective);\\
(a3) $j^A_I : \cI \to SA(\cI)$ satisfies the universal property:  for each  SUR-algebra  $\cS$  and each  pSUR-morphism $h: \cI \to \cS$, there is a unique SUR-morphism $h_A: SA(\cI) \to \cS$  such that $h_A \circ j_I^A = h$. Moreover:\\
\bu If $t'$ is injective, then  $h$ is a $\Sigma-quasi$-embedding iff $h_A$ is a $\Sigma-quasi$-embedding.

(b)  \\
(b1) $j^T_I$ is a fpSUR-morphism and a $\Sigma-quasi$-embedding, that is a $\Sigma$-embedding whenever $\cI$ is transitive;\\
(b2) If $t : C^{t}_s(I) \to I$ is injective (respect. surjective/bijective), then $t^\infty : C^{t^\infty}_s(ST(I)) \to ST(I)$ is injective (respect. surjective/bijective);\\
(b3) $j^T_I : \cI \to ST(\cI)$ satisfies the universal property:  for each   {\bf transitive} SUR-algebra  $\cS$  and each  pSUR-morphism $h: \cI \to \cS$, there is a unique SUR-morphism $h_T: ST(\cI) \to \cS$ such that $h_T \circ j_I^T = h$. 




\ete 

\bdm 

Item (a): based on based on Lemma \ref{IpSUR-le} and  Proposition \ref{fpSURcolim-pr}, 
we can define, by transfinite recursion a convenient {\em increasing} (compatible) family of diagrams $D_\alpha : [0, \alpha] \to pSUR-alg$, $\alpha \in On$, where:\\
(D0) $D_0(\{0\}) = I$;\\
(D1) For each $0 \leq \gamma < \beta < \alpha$, $D_\alpha(\gamma,\beta) = D_\beta(\gamma,\beta) : D_\beta(\gamma) \to D_\beta(\beta)$ is $\Sigma$-embedding and a fpSUR-morphism;\\

Just  define   $D_\alpha(\alpha) = (D_\alpha^{(\alpha)})^+$,  where $D_\alpha^{(\alpha)} := colim_{\beta < \alpha} D_\alpha(\beta)$ and take, for $\beta < \alpha$,  $D_\alpha(\beta,\alpha) = (h^\alpha_\beta)^+ :  D_\alpha(\beta) \to  (colim_{\beta < \alpha} D_\alpha(\beta))^+$ be the unique pSUR-morphism --that is automatically a fpSUR-morphism and a $\Sigma$-embedding, whenever $h^\alpha_\beta$ satisfies this conditions (see Lemma \ref{IpSUR-le}.(d))-- such that $(h^\alpha_\beta)^+ \circ i_0 = h^\alpha_\beta$, where $i_0 : (colim_{\beta < \alpha} D_\alpha(\beta)) \to (colim_{\beta < \alpha} D_\alpha(\beta))^+$ and where $h^\alpha_\beta : D_\alpha(\beta) \to  (colim_{\beta < \alpha} D_\alpha(\beta))$ is the colimit co-cone arrow: by the recursive construction and by Proposition \ref{fpSURcolim-pr} $h^\alpha_\beta$  is a fpSUR-morphism and a $\Sigma$-embedding. This completes the recursion.

Gluing this increasing family of diagrams we obtain a diagram $D_\infty : On \to pSUR-alg$.

By simplicity we will just denote:\\
\bu  $SA(I)_\alpha := D_\infty(\alpha)$, $\alpha \in On$;\\
\bu  $SA(I)_\infty := colim_{\alpha \in On}  SA(I)_\alpha$; \\
\bu $D_\alpha(\beta,\alpha) = j^A_{\beta,\alpha}$, for each $0 \leq \beta \leq \alpha \leq \infty$ (since the family $(D_\alpha)_\alpha$ is increasing, we just have introduce notation for "new arrows").

Then we set: $SA(I) := SA(I)_\infty$ and $j^A_I := j^A_{0,\infty}$.

The verification that $SA(I)$ is a SUR-algebra that satisfies the property in item (a2) and that $j^A_I$ satisfies item (a1)\footnote{In fact, $j^A_{\beta,\alpha}$ is $\Sigma$-embedding whenever $0 \leq \beta \leq \alpha \leq \infty$ and $j^A_{\beta,\alpha}$ is a fpSUR-morphism whenever $0 \leq \beta < \alpha \leq \infty$.},  follows the recursive construction of the diagram and  from a combination of Proposition \ref{fpSURcolim-pr} and Lemma \ref{IpSUR-le}. 

By the same Lemma and Proposition combined, it can be checked by induction that for each $\alpha \in On$, there is a unique  pSUR-morphism $h_\alpha: SA(\cI)_\alpha \to \cS$  such that $h_\alpha \circ j^A_{0,\alpha} = h$ and such that  $h_\alpha$ is injective (respect. $\Sigma-quasi$-embedding, $\Sigma$-embedding), whenever $h$ is  injective (respect. $\Sigma-quasi$-embedding, $\Sigma$-embedding). By applying one more time the colimit construction, we can guarantee that there is  a  unique pSUR-morphism $h^A := h_\infty: SA(\cI)_\infty \to \cS$ such that $h^A \circ j_I^A = h$ and that it satisfies the additional conditions.

\

The proof of item (b) is analogous to the proof of item (a): basically we just have to replace to use of technical Lemma \ref{IpSUR-le} by other technical Lemma  \ref{IpSURtr-le}. In general,  we can on guarantee that $j^T_{\beta,\alpha}$ is a $\Sigma$-embedding and a fpSUR-morphism only for $0 < \beta < \alpha \leq \infty$.
\qdr\edm 

\

In particular, taking $I = \cS_0$ as the {\em initial object} in pSUR-alg (see Proposition \ref{inipSUR-pr} in  Subsection 2.3.3), we have that $SA \cong SA(I)$ and $ST \cong ST(I)$, and they  satisfy  corresponding universal properties:

\bco \label{SASTinitial-co}

(a)  $SA$ is universal (= initial object) over all SUR-algebras,  i.e. for each SUR-algebra $\cS$, there is a unique SUR-algebra morphism $f_S : SA \to \cS$.

(b)  $ST$ is universal (= initial object) over all {\bf transitive} SUR-algebras, i.e. for each {\bf transitive} SUR-algebra $\cS'$, there is a unique SUR-algebra morphism $h_{S'} : ST \to \cS'$.

\eco

\bdm Item (a): Since for each each SUR-algebra $\cS$ there is a unique pSUR-morphism $u_S : \cS_0 \to \cS$ then, by Theorem \ref{adj-te}.(a) above,  $SA(\cS_0)$ is a SUR-algebra that has the required universal property, thus we only have to guarantee that $SA \cong SA(\cS_0)$.  Taking into account the Remark \ref{NOM-re} and the constructions performed in the proof of the item (a) in Theorem above, that we have a (lage) family of compatible pSUR-isomorphisms  $ST_\alpha \cong SA(I)_\alpha, \all \alpha \in On$. Thus 
$SA = \bigcup_{\alpha \in On} SA_\alpha \cong colim_{\alpha \in On} SA(\cS_0)_\alpha = SA(\cS_0)_\infty = SA(\cS_0)$. 

For item (b) the reasoning is similar: note that $I =  \cS_0 = \{*\}$ is a  transitive partial SUR-algebra to conclude that $ST(\cS_0)$ has the required universal property and  note that by the proof of item (b) in Theorem \ref{adj-te} above, that $ST = \bigcup_{\alpha \in On} ST_\alpha \cong colim_{\alpha \in On} ST(\cS_0)_\alpha = ST(\cS_0)_\infty = ST(\cS_0)$.
\qdr\edm

\

			This Corollary describes, in particular, that $SA$ and $ST$ are "rigid" as $\Sigma$-structures and :\\
\bu $SA$ and $C_s(SA)$ are isomorphic SUR-algebras and the universal map $SA \to C_s(SA)$ is the unique iso from $SA$ to $C_s(SA)$;\\
\bu  $ST$ and $C_s(ST)$ are isomorphic SUR-algebras and the universal map $ST \to C_s(ST)$ is the unique iso from $ST$ to $C_s(ST)$.

\

We finish this Section with an application of the Corollary \ref{SASTinitial-co} above: we obtain some non-existence results.

\bco \label{nonexistence-co}\\
(i) Let ${\cal L}$ be a linear SUR-algebra, i.e., $<$ is a total relation (for instance take $\cL =  No$). Then there is no SUR-algebra morphism $h : \cL \to ST$.\\
(ii) Let ${\cal T}$ be a transitive SUR-algebra, i.e., $<$ is a transitive relation (for instance take $\cT = ST, No$). Then there is no SUR-algebra morphism $h : \cT \to SA$.
\eco
\bdm 
(i) Suppose that there is a SUR-algebra morphism $h : \cL \to ST$. Since the binary relation $<$ in $L$ is acyclic and total, it is a strictly linear order, in particular it is transitive. Let $a, b \in L$, since $L$ is linear, $a<b$ in $L \ \Lra \ h(a)<h(b)$ in $ST$.   Now, by the universal property of $ST$ (see Theorem above) there is a unique SUR-algebra morphism $u : ST \to \cL$ and then $h \circ u = id_{ST}$. Summing up, $h : (L,<) \to (ST,<)$ is an isomorphism of structures, thus $(ST,<)$ is a strictly ordered class, but the members of $ST$ $0$ and  $\ll \{-1 \}, \{ 1 \} \rr$ are not comparable  by Fact 1.(c) in the Subsection \ref{SURST-subsec}, a contradiction.

(ii) Suppose that there is a SUR-algebra morphism $h : \cT \to SA$. Since the binary relation $<$ in $T$ is transitive, by the universal property of $ST$ there is a (unique) SUR-algebra morphism $v : ST \to \cT$, thus we get a SUR-algebra morphism $g = h \circ v : ST \to SA$. By the universal property of $SA$ there is a unique SUR-algebra morphism $u : SA \to ST$ (u is a inclusion) and then $g \circ u = id_{SA}$. Thus, for each $a, b \in SA$,  $a<b$ in $SA \ \Lra \ u(a)<u(b)$ in $ST$,  but the members of $SA$ denoted by $-1$ and  $1$ are not related (see Fact.(c) in the Subsection \ref{SURSA-subsec}) and $u(-1) < u(1)$ in $ST$ (by Fact 1.(c) in the Subsection \ref{SURST-subsec}), a contradiction.   \qdr\edm

\section{Conclusion and future works} \label{conclusion}




	


The present work is essentially a collection of elementary  results where we develop, from scratch, a new (we hope!) and  complementary aspect of the Surreal Number Theory. 

In a continuation of the present work (\cite{RM19})
we will establish links, in both directions, 
between SUR-algebras and ZF-algebras (the keystone of  Algebraic Set Theory) and develop the first steps of  a certain kind of set theory based (or ranked)  on surreal numbers, that expands the relation between $V$ and $On$. There is much work to be done: it is clear for us that we just gave the first steps in the Surreal Algebras Theory and in  Set theory based on Surreal Algebras.

In the sequel, we briefly present a (non-exhaustive) list of questions that have occurred to us during the elaboration of this work that we intend to address in the future.

\

{\bf Questions directly connect with the material presented in this work:}

We have described some general constructions in categories of partial SUR-algebra (with at least 2 kinds of morphisms): initial object, non-empty products, substructures and some kinds of directed inductive (co)limits. There are other general constructions available in these categories like quotients and coproducts? A preliminary analysis was made and indicates that the characterizations of the conditions where such constructions exists is a non trivial task.

A specific construction like the (functor) cut surreal for SUR-algebras and its partial version turns out to be very useful to the development of the results of the (partial) SUR-algebra theory: the situation is, in some sense, parallel to the specific construction of rings of fractions construction in Commutative Algebra and Algebraic Geometry. There are other natural and nice specific constructions of (partial) SUR-algebra that, at least, provide new classes of examples?

We have provided, by categorial methods, some universal results that characterizes the SUR-algebras $SA$ and $ST$, and also some relative versions with base ("urelements") $SA(I)$, $ST(I')$ where $I$, $I'$ are partial SUR-algebra satisfying a few constraints. There is an analog result satisfied by the SUR-algebra $No$? There are some natural expansions of $No$ by convenient $I''$ are partial SUR-algebra, $No(I'')$, that also satisfies a universal property that characterizes it up to a unique isomorphism?

\

We saw  that the free/initial  SUR-algebra SA is, in many senses, an expansion of the free/initial ZF-algebra V ans its underlying set theory. Relatively constructions are available for SUR-algebras and for ZF-algebras (see \cite{JMbk}). In particular, it can be interesting examine possible natural expansions of set theories:\\
  (i) with  urelements $B$, $V(B)$,  to some convenient relatively free SUR-algebra $SA(\hat{B})$; \\
	(ii) obtained from the free transitive SUR-algebra $ST \to No$

A combination of the tree lines of research above mentioned can be a interesting ("second-order") task: it will be a line of development of general relative set theories that are base independent.

\

{\bf Unexplored possibilities:}

There exists at least two major aspects of the theory of SUR-algebras that we have not addressed in this work:\\
\bu  the analysis of its model-theoretic aspects;\\
\bu  the consideration of possible applications of SUR-algebras into "traditional" set/class theory,  to answer specific  questions on ZFC/NBG theories.

It is worthy to note that two lines of research can present interesting  cross feedings, as the considerations below will indicate.

First of all, we recall that:\\
-  rational number line $(\mathbb{Q},<)$ is a (or "the") countable dense totally ordered set without endpoints;\\
-  a dense totally ordered set without endpoints is a $\eta_\alpha$-set if and only if it is $\aleph_\alpha$-saturated structure, $\alpha \in On$; \\
- the the surreal number line, $(No,<)$, is for proper class linear orders what the rational number line $(\mathbb{Q},<)$ is for the countable
linear orders.  In fact, $(No,<)$ is a proper class Fra\"iss\'e  limit of the class of all finite linear orders. The surreal numbers are set-homogeneous and universal for all proper class linear orders.\\
- the relational structure $(S,<)$ underlying a SUR-algebra ${\cal S}$  is acyclic and $\eta_\infty$, a natural generalization of the properties above mentioned.

We consider below two remarkable instances of model-theoretic properties applied to set theory that, we believe, could be related to our setting:

{\bf (I)} J. Hamkins have defined in \cite{Ham13}  the notion of "hypnagogic digraph", $(Hg, \rightharpoonup)$, an acyclic digraph graded on $(No,<)$\footnote{I.e., it is given a "rank" function $v : Hg \to No$ such that: for each $x,y \in Hg$, if $x \rightharpoonup y$, then $v(x)< v(y)$.}. 
The hypnagogic digraph is a proper-class analogue the countable random $\mathbb{Q}$-graded digraph: it  is the Fra\"iss\'e limit of the class of all finite $No$-graded digraphs. It is simply the $On$-saturated $No$-graded class digraph,  making it set-homogeneous and universal for all class acyclic digraphs. \\
Hamkins have applied this structure, and some relativized versions,  to prove interesting results concerning models of ZF set theory. For instance:\\
\bu every countable model of set theory $(M, \in^M)$, is isomorphic to a submodel of its own constructible universe $(L^M, \in^M)$;\\
\bu the class of countable models of ZFC is linearly pre-ordered by  the elementary embedding relation.

As a part of a program of model theoretic studies of (relatively free) SUR-algebras, seems natural to determine (and explore) a precise relation between the $No$-ranked relational classes $(Hg, \rightharpoonup)$ and $(SA,<)$ (or $(ST,<)$). And what about the relativized versions of $Hg$ and $SA$ (or $ST$)? This kind of question is very natural as part of an interesting general investigating program relating Model Theory and Category Theory: in one hand we have the model-theoretic universality (from inside or above) of $Hg$ and, on the other hand, we have the category-theoretic universality of (relatively) free constructions (to outside or below) of $SA$ and $ST$.


Can we construct new models of ZF(C) by establishing relations \\
$[Cat^+(SA)] \ \rightsquigarrow \ [Hamkins\ digraph\ models]$ (and/or some variants)\\
in  a way in some sense analogous to the relation:\\
$[sheaves \ over \ boolean \ algebras] \  \rightsquigarrow \ [Cohen \ forcing\ models]$?

{\bf (II)} J. Hirschfeld have provided in \cite{Hir75} a list of axioms - that include axioms for $\in$-acyclicness and for $\in$-density - that describes the model companion of ZF set theory. He emphasizes in the page 369:\\
"{\em \small...This model companion, however, resembles more a theory of order (Theorem 3) than a set theory, and therefore, while supplying an interesting example for model theory it does not shed any new light on set theory...}"

We can wonder about the possible relations of our SUR-set theories and model theoretic (Robinson) forcing(s). This is a natural question since the models of the model companions of ZF have a  "nice" relational structure and the model theoretic forcing can provide the description of model companion/completion of a first order theory.  Considerations involving large infinitary languages  are  also been in sight, since SUR-algebras are $\eta_\infty$ acyclic relational classes.

\end{document}